\documentclass{article}
\usepackage{amsfonts,amsmath,amssymb,amsthm} 
\usepackage{bm}
\usepackage{url}
\usepackage{tikz}
\usetikzlibrary{tikzmark, arrows.meta, calc, positioning}
\usepackage{appendix}
\usetikzlibrary{arrows.meta}

\newcommand{\alert}[2][magenta]%
{{\color{#1}\mbox{[\hspace{-0.4ex}[}#2\mbox{]\hspace{-0.4ex}]}}}

\title{Generalized Hamming weights of additive codes and geometric counterparts}
\author{Jozefien D'haeseleer\thanks{Department of Mathematics: Analysis, Logic and Discrete Mathematics, Ghent University, 9000 Ghent, Belgium. 
E-mail: \texttt{jozefien.dhaeseleer@ugent.be}}  \,\,and\, Sascha Kurz\thanks{Department of Mathematics, University of Bayreuth, 95440 Bayreuth, Germany.}}
\date{}

\def\cB{\mathcal{B}}
\def\cF{\mathcal{F}}
\def\cH{\mathcal{H}}
\def\cS{\mathcal{S}}
\def\cL{\mathcal{L}}
\def\cM{\mathcal{M}}
\def\cP{\mathcal{P}}
\def\FF{\mathbb{F}}
\def\N{\mathbb{N}}

\def\PG{\mathrm{PG}}
\def\supp{\operatorname{supp}}
\def\wt{\operatorname{wt}}
\newcommand{\qbin}[3]{\genfrac{[}{]}{0pt}{}{#1}{#2}_{#3}}

\theoremstyle{definition}
\newtheorem{definition}{Definition}[section]
\newtheorem{remark}[definition]{Remark}
\newtheorem{example}[definition]{Example}

\theoremstyle{plain}
\newtheorem{theorem}[definition]{Theorem}
\newtheorem{lemma}[definition]{Lemma}

\newtheorem{corollary}[definition]{Corollary}

\newtheorem{proposition}[definition]{Proposition}
\begin{document}

\maketitle

\begin{abstract}
    We consider the geometric problem of determining the maximum number $n_q(r,h,f;s)$ of $(h-1)$-spaces in the projective space $\PG(r-1,q)$ such that each subspace of codimension $f$ contains at most $s$ elements. 
    In terms of coding theory, this corresponds to additive codes with a large $f$th generalized Hamming weight. We also consider the dual problem. Here, we determine the minimum number $b_q(r,h,f;s)$ of $(h-1)$-spaces in $\PG(r-1,q)$ such that each subspace of codimension $f$ contains at least $s$ elements. 
    We fully determine $b_2(5,2,2;s)$ as a function of $s$. We additionally give bounds and constructions for other parameters. {For the computational result we partially use extensive integer linear programming computations.}

  \medskip
  
  
  \smallskip
  
\end{abstract}

\section{Introduction}
\label{sec_intro}
It is well known that a linear $[n,k,d]_q$ code $C$ corresponds to a multiset of $n$ points in the projective space $\PG(k-1,q)$ such that each hyperplane contains at most $n-d$ elements. Therefore, instead 
of asking for linear codes with a large minimum Hamming distance $d$ we can also ask for large multisets of points where not too many elements are contained in a hyperplane. If we replace points by 
$(h-1)$-spaces the coding theoretic equivalent is given by additive codes over $\FF_{q^h}$, which are linear over $\FF_q$. Considering multisets of points such that at most $s$ are contained in any subspace 
of codimension $f$ corresponds to linear codes with a large $f$th generalized Hamming weight. Here we want to consider the maximum number $n_q(r,h,f;s)$ of $(h-1)$-subspaces in $\PG(r-1,q)$ with the property 
each subspace of codimension $f$ contains at most $s$ elements. In coding theory terms we are dealing with additive codes that have a large $f$th generalized Hamming weight. 
The special cases where $h$ or $f$ equals $1$ have been extensively studied. 
Outside of this regime not much seems to be known. By taking the complement we can relate our problem to a 
dual problem: what is the minimum number $b_q(r,h,f;s)$ of $(h-1)$-spaces in $\PG(r-1,q)$ where each subspace of codimension $f$ contains at least $s$ elements? If $s=1$ one also speaks of blocking 
sets of $(h-1)$-spaces w.r.t.\ $(r-f-1)$-spaces. The case $h=f=1$ is a classical problem, and we e.g.\ have $b_q(3,1,1;1)=q+1$ attained by all points on a line in $\PG(2,q)$, which is also called a trivial blocking set. For non-trivial blocking sets, which are those not containing a full line in its support, the minimum size rises to $3(p+1)/2$ for odd primes $p$ \cite{blokhuis1994size}.  For $h=2$ and $r-f=3$ we 
refer to e.g.\ \cite{eisfeld1997blocking,kovacs2025blocking,metsch2004blocking}. If $s>1$ one speaks of multiple blocking sets, see e.g.\ \cite{bishnoi2017minimal}. In \cite[Definition 2.1]{chen2025t} the authors speak of an $s$-fold $f$-blocking set of $\PG(r-1,q)$ for the special case $h=1$. 
In this paper, we use the following definition. 
\begin{definition}
  A multiset $\cM$ of $(h-1)$-spaces in $\PG(r-1,q)$ is called an $s$-fold blocking set w.r.t.\ $(v-1)$ spaces if every $(v-1)$-space in $\PG(r-1,q)$ contains at least $s$ elements from $\cM$.
\end{definition}
Whenever the parameters are clear from the context, we just speak of generalized blocking sets.

The remaining part of this paper is structured as follows. In Section~\ref{sec_preliminaries} we introduce the necessary preliminaries. The relation between the geometric objects and coding theory is outlined in Section~\ref{sec_coding_theory}. In Section~\ref{sec_asymptotic} we summarize our knowledge on the asymptotic behavior of $n_q(r,h,f;s)$. General constructions are studied in Section~\ref{sec_general_constructions}. In Section \ref{sec_blocking_sets} we investigate the generalized blocking sets and their minimum possible size $b_q(r,h,f;s)$. 
 In Section~\ref{sec_n_q_4_2_2_s} we study the maximum number of lines in $\PG(4,q)$ such that each plane contains at most $s$ lines, and we fully determine $b_2(5,2,2;s)$ as a function of $s$. We close with a conclusion and a few open problems in Section~\ref{sec_conclusion}. Sporadic blocking sets, found by integer linear programming searches, are listed in Appendix~\ref{sec_bs_lists}.

\section{Preliminaries}
\label{sec_preliminaries}
The set of all subspaces of $\FF_q^r$, ordered by the incidence relation
$\subseteq$, is called the $(r-1)$-dimensional projective geometry over
$\FF_q$ and denoted by $\PG(r-1,q)$. Here we use the projective dimension, so that an $(i-1)$-space in $\PG(r-1,q)$ 
is an $i$-dimensional space in the vector space setting 
$\FF_q^r$. We will call $0$-, $1$-, $2$-, $3$-, 
and $(r-2)$-spaces points, lines, planes, solids, and hyperplanes, respectively. For two subspaces $S$ and $S'$ we write $S\subseteq S'$ if $S$ is contained in $S'$. Moreover,
we say that $S$ and $S'$ are incident if and only if $S\subseteq S'$ or $S\supseteq S'$.
Let $[i]_q:=\tfrac{q^i-1}{q-1}$ denote the number of points
of an arbitrary $(i-1)$-space in $\PG(r-1,q)$ where $r\ge i$. By convention we
set $[0]_q:=0$. More generally, by $\qbin{r}{i}{q}:=\tfrac{\prod_{j=0}^{i-1} q^{r-j}-1}{\prod_{j=0}^{i-1} q^{i-j}-1}=\tfrac{\prod_{j=0}^{i-1} [r-j]_q}{\prod_{j=0}^{i-1} [i-j]_q}$ we denote the number of $(i-1)$-spaces in $\PG(r-1,q)$. Duality implies $\qbin{r}{i}{q}=\qbin{r}{r-i}{q}$.
It is known that the number of $j$-spaces disjoint from a fixed $m$-space in $\PG(n, q)$ equals $q^{(m+1)(j+1)}\qbin{n-m}{j+1}{q}$, see \cite[Section 170]{Segre}. 

We can represent an $(i-1)$-space in $\PG(r-1,q)$ by an $i\times r$ generator matrix over $\FF_q$. A multiset of points $\cM$ in $\PG(r-1,q)$ is a mapping from the set of points to $\mathbb{N}$. For a given point $P$ we call $\cM(P)$ its multiplicity. We say that $\cM$ is spanning if the points with positive multiplicity span the entire ambient space. The notion of the point multiplicities is extended additively to any subspace $S$ via $\cM(S):=\sum_{P\in S} \cM(P)$. The relation between multisets of points and linear codes is explained in detail in Section~\ref{sec_coding_theory}. For additive codes we need the following generalization, see e.g.\ \cite{ball2024additive}.  

\begin{definition}
  A projective $h-(n, r, s)_q$ system is a multiset $\cS$ of $n$ subspaces of $\PG(r-1,q)$ of dimension at most $(h-1)$ such that each hyperplane contains at most $s$ elements of $\cS$, and some hyperplane contains exactly $s$ elements of $\cS$. We say that $\cS$ is faithful if all its elements have dimension $(h-1)$. A projective $h-(n, r, s)_q$ system $\cS$ is a projective $h-(n, r, s, \mu)_q$ system if each point is contained in at most $\mu$ elements from $\cS$, and there is some point that is contained in exactly $\mu$ elements from $\cS$.
\end{definition}

A faithful projective $1-(n,r,s)_q$ system $\cS$ is just a multiset of points with cardinality $n$ in $\PG(r-1,q)$ with the property that the maximum hyperplane multiplicity $\cS(H)$ equals $s$. Unfaithful projective $h-(n,r,s)_q$ systems also allow the containment of $(-1)$-dimensional subspaces, which correspond to zero columns in the generator matrix of a corresponding linear code for $h=1$, see Section~\ref{sec_coding_theory}. 

\begin{definition}
  By $n_q(r,h;s)$ we denote the maximum number $n$ such that a projective $h-(n, r, s)_q$ system exists.
\end{definition}

Note that the elements of $\cS$ span the entire ambient space $\PG(r-1, q)$ if and only if $s<n$. If $\cS$ is a projective $h-(n, r, s)_q$ system that is not faithful, then we can easily construct a faithful projective $h-(n, r, \le s)_q$ system $\cS'$ by replacing each element $S\in\cS$ with dimension smaller than $h-1$ by an arbitrary $(h-1)$-space containing $S$. The functions $n_q(r,h;s)$ were e.g.\ studied in \cite{kurz2024additive}, and indirectly in any paper on additive codes with good parameters. 

For our situation we need an even more general notion.

\begin{definition}
  \label{def_system_gen}
  A projective $(h,f)-(n,r,s)_q$ system, where $h+f\le r$, is a multiset
  $\cS$ of $n$ subspaces of $\PG(r-1, q)$ of dimension at most $(h-1)$ such
  that each subspace of codimension $f$ contains at most $s$ elements
  of $\cS$, and some subspace of codimension $f$ contains exactly $s$
  elements of $\cS$. We say that $\cS$ is faithful if all elements have
  dimension $(h-1)$. A projective $(h,f)-(n,r,s)_q$ system $\cS$ is a projective
  $(h,f)-(n,r,s,\mu)_q$ system if each $(f-1)$-space is contained
  in at most $\mu$ elements from $\cS$, and there is some $(f-1)$-space that is
  contained in exactly  $\mu$ elements from $\cS$.
\end{definition}

A projective $(h,1)-(n,r,s)_q$ system is just a projective $h-(n,r,s)_q$ system, and a general projective $(h,f)-(n,r,s)_q$ system corresponds to an additive $[n,r/h,d_f]_q^h$ code $C$ with $s=n-d_f$. Here $d_f$ denotes the minimum  $f$th generalized Hamming weight of $C$, see Section~\ref{sec_coding_theory}.
The parameter $\mu$ corresponds to the maximum column multiplicity of linear codes over $\FF_q$, if we identify linear dependent non-zero columns of a generator matrix. 

\begin{definition}
  By $n_q(r,h,f;s)$ we denote the maximum number $n$ such that a projective $(h,f)-(n, r, s)_q$ system exists.
\end{definition}


Clearly we can convert any given projective $(h,f)-(n,r,s)_q$ system into a faithful projective $(h,f)-(n,r,\le s)_q$ system by replacing each element $U$ by an arbitrary $(h-1)$-space containing $U$.

\bigskip 

By $\dim(U)$ we denote 
the projective dimension 
of a subspace $U$ in $\FF_q^n$, which is one less than the algebraic dimension. With this, the subspace distance is given by $d_S(U,V)=\left(\dim(U)+1\right)+\left(\dim(V)+1\right)-2(\dim(U\cap V)+1)=\dim(U)+\dim(V)-2\dim(U\cap V)$, {which is an even number if $\dim(U)=\dim(V)$.} 

By $A_q(v,k;2\delta)$ we denote the maximum number of $(k-1)$-spaces in $\PG(v-1,q)$ with minimum subspace distance $2\delta$. 
Here we have $\dim(U\cap V)+1\le k-\delta$ and speak of constant-dimension codes. In the following we give an upper bound  and one simple construction for constant-dimension codes. In order to keep the paper self-contained we give a brief proof and description. For more details we refer to the survey \cite{kurz2021constructions}. 
\begin{lemma}
  \label{lemma_anticode_bound}
  For $v\ge 2k$ we have
  $$
    A_q(v,k;2\delta)\le \frac{\qbin{v}{k-\delta+1}{q}}{\qbin{k}{k-\delta+1}{q}}.
  $$
\end{lemma}
\begin{proof}
  Since $\dim(U\cap V)+1\le k-\delta$ for any two different elements $U$ and $V$ of the constant-dimension code, each $(k-\delta)$-space is contained in at most one $(k-1)$-space from the constant-dimension code.
\end{proof}
A rank metric code $M$ is a subset of $m\times n$ matrices over $\FF_q$ equipped with the rank distance  $d_r(M,M')=\operatorname{rk}(M-M')$. Assuming $m\le n$, a Singleton-like upper bound is known and gives $|M|\le q^{n(m-\delta+1)}$ for minimum rank distance $\delta$ \cite{delsarte1978bilinear}. Codes attaining this bound are called maximum rank distance (MRD) codes. They exist for all parameters, even if one additionally assumes that the matrices form a linear space, i.e.\ assuming that the code is linearly closed. For a survey on MRD codes we refer to \cite{sheekey201913}. Given an MRD code $M$ of $m\times n$ matrices over $\FF_q$ with minimum rank distance $\delta$, we obtain a lifted MRD (LMRD) code $\cM$ by prepending $m\times m$ unit matrices. Interpreted as generator matrices of $(m-1)$-spaces in $\PG(n+m-1,q)$, $\cM$ is a set of $q^{n(m-\delta+1)}$ $(m-1)$-spaces in $\PG(n+m-1,q)$ such that the dimension of the intersection of any two elements is at most $(m-\delta-1)$ and there exists a special $(n-1)$-space 
$S$ that is disjoint to all elements of $\cM$.  
 
\section{Relation to coding theory}
\label{sec_coding_theory}
A linear $[n,k]_q$ code $C$ is a $k$-dimensional subspace of the vector space $\FF_q^n$. The elements of $C$ are called codewords and the Hamming weight $\wt(c)$ of a codeword $c\in C$ is the number of non-zero entries. With this, the Hamming distance $d\!\left(c_1,c_2\right)$ between two codewords is given by $\wt\!\left(c_1-c_2\right)$. The minimum Hamming distance $d(C)$ of a (linear) code is the minimum Hamming distance $d\!\left(c_1,c_2\right)$ between two different codewords. We say that an $[n,k]_q$ code $C$ is an $[n,k,d]_q$ code if its minimum Hamming distance $d(C)$ equals $d$. A linear code is called $\Delta$-divisible if the weights of all codewords are divisible by $\Delta$. If the non-zero weights of a linear $[n,k]_q$ code are contained in $\left\{w_1,\dots,w_l\right\}$ we also speak of an $\left[n,k,\left\{w_1,\dots,w_l\right\}\right]_q$ code, and an $l$-weight code if all weights are attained. 
As a representation for a linear code we use a $k\times n$ generator matrix over $\FF_q$. The dual code $C^\perp$ of an $[n,k]_q$ code $C$ is the $[n,n-k]_q$ code whose codewords are orthogonal to all codewords in $C$. By $d^\perp$ we denote the corresponding minimum Hamming distance. We say that $C$ has full length if $d^\perp\ge 2$, which is equivalent to the property that there is no zero-column in a given generator matrix for $C$. 
It is well known that full length $[n,k]_q$ codes are in one-to-one correspondence to spanning multisets of cardinality $n$ in $\PG(k-1,q)$, see e.g.\ \cite{dodunekov1998codes}.\footnote{ Given a linear $[n,k]_q$ code $C$ with generator matrix $G$, we can interpret its columns as $1$-dimensional vector spaces of $\FF_q^k$ or points in $\PG(k-1,q)$.}
The minimum Hamming distance $d$ of a linear code corresponds to the geometric property that the maximum number of elements of the multiset of points that is contained in a hyperplane is given by $n-d$. So, a large minimum Hamming distance corresponds to a small maximum number of points in hyperplanes. Alternatively, minimizing the possible length $n$ of an $[n,k,d]_q$ is equivalent to maximizing the cardinality of a multiset of points in $\PG(k-1,q)$ with at most $s$ points in each hyperplane, where $s=n-d$. 

More generally, a block code $C$ of length $n$ over the alphabet $\FF_q$ is just a subset of $\FF_q^n$ (equipped with the Hamming metric). If $C$ is linearly closed, i.e.\ if $c,c'\in C$ and $\alpha,\beta\in\FF_q$ implies $\alpha c+\beta c'\in C$, then we have a linear $[n,k]_q$ code, where $k=\log_q |C|$ is called the dimension. An additive code is just a block code that is additively closed, i.e.\ $c,c'\in C$ implies $c+c'\in C$. Each additive code is linear over some subfield, see e.g.\ \cite{ball2023additive}. By an $[n,r/h,d]_q^h$ code we denote an additive code $C\subseteq \FF_{q^h}^n$ that is linear over $\FF_q$, has minimum Hamming distance $d$ and cardinality $q^r$. We call $r/h\in\mathbb{Q}$ its dimension. We can represent an $[n,r/h,d]_q^h$ code as the $\FF_q$ row span of an $r\times n$ generator matrix $G$ over $\FF_{q^h}$. Choosing an $\FF_q$ basis of $\FF_{q^h}$ we can expand this generator matrix to a subfield generator matrix $\widetilde{G}\in\FF_q^{r\times nh}$. By $\mathcal{X}_G(C)$ we define the multiset of the $n$ subspaces spanned by the
$n$ blocks of $h$ columns of $\widetilde{G}$ in this way. Note that these subspaces give rise to projective systems.

\begin{theorem}(\cite[Theorem 5]{ball2024additive})
   \label{thm_connection}
   If $C$ is an additive $[n,r/h,d]_q^h$ code with generator matrix $G$, then
   $\mathcal{X}_G(C)$ is a projective $h-(n, r, n-d)_q$ system $\mathcal{S}$, and conversely,  each projective $h-(n,r, s)_q$ system $\mathcal{S}$ defines an additive $[n,r/h,n-s]_q^h$ code $C$.
\end{theorem}

The parameters of a linear $[n,k,d]_q$ code $C$ are related by the so-called
\emph{Griesmer bound} \cite{griesmer1960bound,solomon1965algebraically}
\begin{equation}
  \label{eq_griesmer_bound}
  n\ge \sum_{i=0}^{k-1} \left\lceil\frac{d}{q^i}\right\rceil=:g_q(k,d).
\end{equation}
From this one can derive the bound 
\begin{eqnarray}
  n&\ge& \left\lceil
  \frac{g_q\!\left(r,d\cdot q^{h-1}\right)}{[h]_q}
  \right\rceil
  =
  \left\lceil
  \frac{ \sum\limits_{i=0}^{r-1} \left\lceil d\cdot q^{h-1-i}\right\rceil}{[h]_q}
  \right\rceil\nonumber\\ 
  &=&d+\left\lceil\frac{\sum\limits_{i=1}^{r-h} \left\lceil\frac{d}{q^i}\right\rceil}{[h]_q}\right\rceil
  =
  d+ \left\lceil \frac{g_q(r-h+1,d)-d   }{[h]_q}\right\rceil
\end{eqnarray}
for the parameters of an additive $[n,r/h,d]_q^h$ code, see e.g.\ \cite[Theorem 12]{ball2024additive} or \cite[Lemma 15]{kurz2024additive}. Using Theorem~\ref{thm_connection} this gives an upper bound for $n_q(r,h,1;s)$, which we call the Griesmer upper bound. More precisely, we call the largest integer $n$ that satisfies $[h]_q\cdot n\ge g_q\!\left(r,(n-s)\cdot q^{h-1}\right)$ the Griesmer upper bound for $n_q(r,h;s)$, see e.g.\ \cite[Example 5]{kurz2024additive}.

The Hamming weight $\wt(c)$ turns $\FF_q^n$ into a normed vector space. For $c=\left(c_1,\dots,c_n\right)\in\FF_q^n$ we call
\begin{equation}
  \supp(c):=\left\{1\le i\le n\,:\, c_i\neq 0\right\}
\end{equation}
the support of $c$, so that $\wt(c)=|\supp(c)|$. 
For some linear subspace $C$ in $\FF_q^n$ let
\begin{equation}
  \supp(C):=\left\{1\le i\le n\,:\, \exists c=\left(c_1,\dots,c_n\right)\in C, c_i\neq 0\right\}
\end{equation}
be the support of $C$ and $\dim_{\FF_q}(C)$ its (algebraic) $\FF_q$-dimension. For two
$\FF_q$ vector spaces $C$, $C'$ in $\FF_q^n$ we write $C\subseteq C'$ if $C$ is contained in $C'$. With this, the $f$th generalized Hamming weight of
a linear code $C$ \cite{helleseth1977weight,klove1978weight}, denoted as $d_f(C)$, is the size of the
smallest support of an $f$-dimensional subcode of $C$:
\begin{equation}\label{eq:dflinear}
   d_f(C):= \min\!\left\{|\supp(C')|\,:\, C'\subseteq C, \dim_{\FF_q}(C')=f\right\}.
\end{equation}
In particular, $d_1(C)$ is the minimum Hamming distance of a linear code $C$. The sequence $\left(d_1(C),\dots,d_k(C)\right)$ is called the weight hierarchy of a linear $[n,k]_q$ code $C$. Clearly, we have $1\le d_1(C)\le \dots\le d_k(C)\le n$. The generalized Hamming weights can be used to describe the cryptographic performance of a linear code over the wire-tap channel of type~II \cite{wei1991generalized}. Moreover, it van also be used to determine the trellis complexity of the code \cite{chen2001trellis,forney1994density,forney1994dimension,kasami1993optimum}. The weight hierarchy of a linear code can be obtained from a quadratic form over a finite field \cite{li2020weight,li2022weight,liu2023generalized}. Also the geometric reformulation of the generalized Hamming weights in terms of multisets of points is well known \cite{helleseth1992generalized,tsfasman1995geometric}. Let $\cM$ be a multiset of points in $\PG(k-1,q)$ and $C$ its corresponding $[n,k]_q$ code. Then, we have
\begin{equation}
  \label{eq_gen_ham_geometry}
  n-d_f(C) =\max\left\{\cM(U)\,:\, U\text{ subspace of codimension }f \right\}
\end{equation} for all $1\le f\le k$. In order to keep the paper self-contained we state a brief argument, c.f.\ \cite{ball2025additive}.  
Given a linear $[n,k,d]_q$ code $C$, a codeword (a $1$-dimensional subcode) of $C$ is obtained by left multiplication of a generator matrix $G$ by a vector $v\in\FF_q^k$. Considering $v$ as a point in $\PG(k-1, q)$, the hyperplane $v^\perp$ 
contains the point $x$ if and only if $\langle v,x\rangle=0$. Therefore, if we take the set of $n$ points in $\PG(k-1,q)$, corresponding to the columns of $G$, we have that the codeword $vG$ has weight $w$ if and only if $n-w$ of these points are contained in the hyperplane $v^\perp$. 
More generally, for a $j$-dimensional subspace $V$ of $\FF_q^n$ the codimension $j$ subspace $V^\perp$ contains $n-w$ points if the subspace $\left\{vG\,:\, v\in V\right\}$ has support size $w$, which proves Equation~(\ref{eq_gen_ham_geometry}). We can apply the same argument to the subfield generator matrix $\widetilde{G}$ of an additive $[n,r/h,d]_q^h$ code $C$ to conclude 
\begin{equation}
  \label{eq_gen_ham_geometry_add}
  n-d_f(C) =\max\left\{\left|\left\{S\in\mathcal{X}_G(C)\,:\, S\le U\right\}\right|\,:\, U\text{ subspace of codim. }f \right\}
\end{equation} for all $1\le f\le k$. Hence, looking for good additive codes, corresponds to look for large  projective systems.
\begin{theorem}(Griesmer-type bound) \cite[Theorem 4]{helleseth1995bounds}, \cite[Theorem 5]{helleseth1992generalized}\\\label{thm_griesmer_gen_ham}
   For each linear $[n,k]_q$ code and each $1\le f\le k$ we have
   \begin{equation}
     \label{ie_griesmer_gen_ham_weight}
     n \ge d_f +\sum_{j=1}^{k-f} \left\lceil \frac{d_f}{[f]_q\cdot q^j}\right\rceil=:g_q^f\!\left(k,d_f\right).
   \end{equation}
\end{theorem}

Currently we do not know any Griesmer type bound for the $f$th generalized Griesmer weight of additive codes, which is tight for all sufficiently large minimum distances. This is an important open problem. For $f>1$ and $h>1$ we cannot reconstruct the number of $(h-1)$-spaces in a codimension $f$ space from the number of contained points and the total number of $(h-1)$-spaces. 

\section{Asymptotic results}
\label{sec_asymptotic}
We first state a sum construction and an easy upper bound that can be asymptotically attained. 
\begin{lemma}
   \label{lemma_sum_gen_ham}
   $n_q\!\left(r,h,f;s_1+s_2\right)\ge n_q\!\left(r,h,f;s_1\right)+n_q\!\left(r,h,f;s_2\right)$ 
\end{lemma}
\begin{proof}
  Consider the union of a projective $(h,f)-\left(n_q(r,h,f;s_1),r,s_1\right)_q$ and a projective $(h,f)-\left(n_q(r,h,f;s_2),r,s_2\right)_q$ system.
\end{proof}

\begin{lemma}
\label{lemma_one_weight_bound_gen_hamming_weight}
We have
\begin{equation}
  n_q(r,h,f;s)\le \frac{\qbin{r}{f}{q}\cdot s}{\qbin{r-h}{f}{q}}
  =
  \prod_{i=0}^{f-1} \frac{[r-i]_q}{[r-h-i]_q}\cdot s.
\end{equation}
\end{lemma}
\begin{proof}
   Let $\cS$ be a faithful projective $(h,f)-(n,r,s)_q$ system with
  $n=n_q(r,h,f;s)$. Since each element $S\in\cS$ is contained in
  $\qbin{r-h}{f}{q}$ subspaces of codimension $f$ and there are $\qbin{r}{f}{q}$ subspaces of codimension $f$ in total,
  we conclude $n\le \tfrac{\qbin{r}{f}{q}\cdot s}{\qbin{r-h}{f}{q}}$.
\end{proof}

Considering the set of all $n=\qbin{r}{h}{q}$ $h$-spaces in $\PG(r-1,q)$ we see that the upper bound in Lemma~\ref{lemma_one_weight_bound_gen_hamming_weight} is tight for $s=\qbin{r-f}{h}{q}$. Using $\lambda$ copies of this construction yields 
\begin{equation}
  \lim_{s\to\infty} n_q(r,h,f;s) \cdot \frac{\qbin{r-h}{f}{q}}{\qbin{r}{f}{q}\cdot s}=1.
\end{equation}
For $f=1$ the Griesmer bound implies that the difference between $n_q(r,h,1;s)$ and the corresponding Griesmer upper bound tends to zero if $s$ tends to infinity, which is a much tighter statement.\footnote{While the geometric equivalent of linear or additive codes is very handy for many situations, here the coding theory version looks more nicely. In particular, denoting the minimum length $n$ of an $[n,k,d]_q$ code by $\tilde{n}_q(k,d)$, we have $\lim_{d\to\infty} \tilde{n}_q(k,d)-g_q(k,d)=0$. There is a similar formulation for additive codes, see \cite{kurz2024additive}.} The same stronger result also holds for the cases where $h=1$ and $f$ is arbitrary. Of course it would be very interesting to have such a result in general. However, the relation to constant-dimension codes in  Lemma~\ref{lemma_connection_to_subspace_codes} indicates that this might be a hard problem.  

Instead of letting $s$ tend to infinity we can also consider $n_q(r,h,f;s)$ as a sequence in the field size $q$. 
\begin{lemma}
  \label{lemma_connection_to_subspace_codes}
  For $1\le \delta\le h$ and $r\ge 2h$ we have 
  $$
    n_q(r,h,r-h-\delta+1;1)=A_q(r,h;2\delta).
  $$  
\end{lemma}
\begin{proof}
Since each $(h+\delta-2)$-space contains at most one $(h-1)$-space from the projective system, the dimension formula implies that each $(h-\delta)$-space is contained in at most one $(h-1)$-space from the projective system. With this, the distance between the $(h-1)$-spaces $U$ and $V$ is $\dim(U)+\dim(V)-2\dim(U\cap V) \geq 2h-2-2(h-\delta-1) = 2\delta$.
\end{proof}
The special case $\delta=h$ corresponds to partial spreads where many bounds are known, see e.g.\ \cite[Section 9]{kurz2021divisible}. For general parameters the following construction using (L)MRD codes is well known.
\begin{proposition}
  \label{prop_lmrd_construction}
  For $1\le \delta\le h$ and $r\ge 2h$ we have 
  $$
    n_q(r,h,r-h-\delta+1;1)\ge q^{(r-h)(h-\delta+1)}. 
  $$
\end{proposition}
\begin{proof}
  Consider an LMRD code $\cM$ of $q^{(r-h)(h-\delta+1)}$ $(h-1)$-spaces in $\PG(r-1,q)$ with minimum rank distance $\delta$, i.e., the projective dimension of the intersection of two different elements of $\cM$ is at most $(h-\delta-1)$. Thus, each subspace of codimension $f=r-h-\delta+1$ contains at most one element from $\cM$.
\end{proof}
As an example, for $r=6$, $h=3$, and $\delta=2$ we obtain a set of $q^6$ planes in $\PG(5,q)$ with the property each solid contains at most one plane (and each $4$-space contains at most $q^3$ planes). Via Lemma~\ref{lemma_connection_to_subspace_codes} we can replace the used LMRD codes by any other constant-dimension code with the same minimum subspace distance, see e.g.\ \cite{kurz2021constructions} for a survey on some constructions from the literature. For the aforementioned parameters we remark that this construction is not optimal, since there is a construction known with $q^6+2q^2+q+1$ planes, see \cite{hkk2015}.

\begin{corollary}
  For $1\le \delta\le h$ and $r\ge 2h$ we have 
  $$
    n_q(r,h,r-h-\delta+1;s)\ge s\cdot q^{(r-h)(h-\delta+1)}. 
  $$
\end{corollary}

For the special cases where either $r=2h$ or $s=1$ the construction with the LMRD codes is asymptotically tight if $q$ tends to infinity.

To prove this, we use the so-called $q$-Pochhammer symbol $$(a;q)_n=\prod_{i=0}^{n-1} \left(1-aq^i\right).$$ In particular, we will use that \begin{equation}
    \label{ie_qbin}
    1\le q^{-b(a-b)}\cdot \qbin{a}{b}{q}
    \le \frac{1}{(1/q;1/q)_b}, 
  \end{equation}
  see e.g.\ \cite{koetter2008coding}.
\begin{proposition}
  \label{prop_asymtotic_q_1}
  For $1\le \delta\le h$ we have
  $$
    \lim_{q\to\infty} \frac{n_q(2h,h,h-\delta+1;s)}{s\cdot q^{h(h-\delta+1)}}=1.
  $$
\end{proposition}
\begin{proof}
Lemma~\ref{lemma_one_weight_bound_gen_hamming_weight} yields 
  $$
    n_q(2h,h,h-\delta+1;s)\le \frac{\qbin{2h}{h-\delta+1}{q}\cdot s}{\qbin{h}{h-\delta+1}{q}}.  
  $$  
  Applying the $q$-Pochhammer symbols to both numerator and denominator, and using \eqref{ie_qbin}, we have
  \begin{eqnarray*}
    n_q(2h,h,h-\delta+1;s)&\le& \frac{q^{(h+\delta-1)(h-\delta+1)}}{q^{(\delta-1)(h-\delta+1)}} \cdot \frac{1}{(1/q;1/q)_{h-\delta+1}}\cdot s\\ 
    &=& 
    q^{h(h-\delta+1)} \cdot \frac{1}{(1/q;1/q)_{h-\delta+1}}\cdot s,
  \end{eqnarray*}
  where the second factor tends to $1$ as $q$ approaches infinity.
\end{proof}

As in the proof of Proposition~\ref{prop_asymtotic_q_1},   Lemmas~\ref{lemma_anticode_bound} and \ref{lemma_connection_to_subspace_codes}, together with Inequality~(\ref{ie_qbin}) yield the upper bound, while Proposition \ref{prop_lmrd_construction} provides a matching lower bound.
\begin{proposition}
  For $1\le \delta\le h$ and $r\ge 2h$ we have 
  $$
    \lim_{q\to\infty} \frac{n_q(r,h,r-h-\delta+1;1)}{q^{(r-h)(h-\delta+1)}}=1. 
  $$
\end{proposition}
We remark that Lemma~\ref{lemma_anticode_bound} is known as the anticode bound in the context of subspace codes and that tighter bounds are known, see e.g.\ \cite{kurz2021constructions}. 

\section{General constructions}
\label{sec_general_constructions}
In this section we want to study known constructions for linear codes from the literature and generalize them to our context.

In coding theory it is well known that the problem of determining the minimum possible length of an $[n,k,d]_q$ code as a function of $d$ is a finite problem for given parameters $k$ and $q$. More precisely, if the minimum distance $d$ is sufficiently large, then the Griesmer bound can always be attained with equality. A corresponding construction was given by Solomon and Stiffler \cite{solomon1965algebraically}. In geometric terms this means that the determination of the function $n_q(r,1,1;\cdot)$ in terms of $s$ is a finite, but still rather hard, problem for each given pair of parameters $r$ and $q$. In \cite{kurz2024additive} this result was generalized to additive codes, i.e.\ also applies to $n_q(r,h,1;\cdot)$ for arbitrary $h$. In order to describe the Solomon--Stiffler construction and its generalization, we have to introduce further notation. For each subspace $S$ in $\PG(r-1,q)$ we denote its characteristic function by $\chi_S$, i.e.\ we have, for a point $P\in PG(r-1,q)$ that $\chi_S(P)=1$ if $P\in S$ and $\chi_S(P)=0$ otherwise.
\begin{definition}
  We say that a multiset of points $\cM$ in $\PG(r-1,q)$ is $h$-partitionable if there exist $(h-1)$-spaces $S_1,\dots,S_l$, for some integer $l$, such that $\cM=\sum_{i=1}^l \chi_{S_i}$, i.e.\ $\cM$ can be partitioned into $(h-1)$-spaces.
\end{definition}
To ease the notation and to avoid technical difficulties, we choose a chain of subspaces $S_1\subseteq S_2\subseteq\dots\subseteq S_r$ in $\PG(r-1,q)$, where $S_i$ has projective dimension $(i-1)$.

\begin{definition}
  \label{def_partitionable}
  Given a chain of subspaces $S_1\subsetneq S_2\subsetneq\dots\subsetneq S_r$ in $\PG(r-1,q)$, we say that $\sum_{i=1}^r a_iS_i$ is $h$-partitionable over $\FF_q$ if the multiset of 
  points $\sum_{i=1}^r a_i\chi_{S_i}$ in $\PG(r-1,q)$ is $h$-partitionable, where $a_i\in\mathbb{Z}$ for all $1\le i\le r$.
\end{definition}
For example, we trivially have that $S_3$ is $3$-partitionable over $\mathbb{F}_q$. The existence of plane spreads in $\mathrm{PG}(5,q)$ implies that $S_6$ is $3$-partitionable over $\mathbb{F}_q$. Since $[7]_q$ is not divisible by $[3]_q$, we have that $S_7$ is not $3$-partitionable over $\mathbb{F}_q$, while $(q^2 + q + 1)\cdot S_7$ is $3$-partitionable.  
For the details on the underlying constructions we refer to \cite{kurz2024additive}.

Using a specific parameterization of the minimum distance $d$ the
Griesmer bound in Inequality~(\ref{eq_griesmer_bound}) can be written more
explicitly as follows.
  Let $k$ and $d$ be positive integers. Write $d$ as
  \begin{equation}
    \label{eq_griesmer_representation_min_dist}
    d=\sigma q^{k-1}-\sum_{i=1}^{k-1}\varepsilon_iq^{i-1},
  \end{equation}
  where $\sigma\in\N_{>0}$, 
  and the $0\le\varepsilon_i<q$ are integers for all $1\le i\le k-1$. Then, Inequality~(\ref{eq_griesmer_bound})
  is satisfied with equality if and only if
  \begin{equation}
    \label{eq_griesmer_representation_length}
    n=\sigma[k]_q-\sum_{i=1}^{k-1}\varepsilon_i[i]_q,
  \end{equation}
  which is equivalent to
  \begin{equation}
    \label{eq_griesmer_representation_species}
    n-d=\sigma[k-1]_q-\sum_{i=1}^{k-1}\varepsilon_i[i-1]_q.
  \end{equation}
\begin{remark}
  Given $k$ and $d$, Equation~(\ref{eq_griesmer_representation_min_dist})
  always determines $\sigma$ and the $\varepsilon_i$ uniquely. This is
  different for Equation~(\ref{eq_griesmer_representation_species}) given
  $k$ and $n-d=s$.


  By relaxing to
  $0\le \varepsilon_i\le q$ we can ensure existence and uniqueness are
  enforced by additionally requiring $\varepsilon_j=0$ for all $j<i$
  where $\varepsilon_i=q$ for some $i$. The same is true for
  Equation~(\ref{eq_griesmer_representation_length}) given $k$ and $n$.
  For more details, we refer to \cite[Chapter 2]{govaerts2003classifications},
  which also gives pointers to Hamada's work on minihypers. We will mostly state our corresponding results referring to Equation~(\ref{eq_griesmer_representation_min_dist}) and using the coding theoretic formulation.
\end{remark}

Given arbitrary $\varepsilon_1,\dots,\varepsilon_{k-1}\in\mathbb{Z}$, we have that
$
  \cM=\sigma \chi_{S_k}-\sum_{i=1}^{k-1}\varepsilon_i\chi_{S_i}
$
is a multiset of points in $\PG(k-1,q)$, which is a projective $1-(n,k,s)_q$ system for all sufficiently large $\sigma\in\N$. Here $n=\sigma [k]_q-\sum_{i=1}^{k-1}\varepsilon_i[i]_q$ and $s=\sigma [k-1]_q-\sum_{i=1}^{k-1}\varepsilon_i[i-1]_q$, see e.g.\ \cite[Lemma 23]{kurz2024additive}. While $\sigma S_k-\sum_{i=1}^{k-1}\varepsilon_i S_i$ is obviously $1$-partitionable over $\FF_q$ if $\sigma$ is sufficiently large, there are further conditions for being $h$-partitionable when $h>1$ 
as well as more sophisticated constructions for the partition, see \cite{kurz2024additive}. 

Recall that a multiset $\cM$ of $(h-1)$-spaces in $\PG(r-1,q)$ is an
$s$-fold blocking set w.r.t.\ $(v-1)$ spaces if every $(v-1)$-space in $\PG(r-1, q)$ contains at least $s$ elements from $\cM$. 
Note that the smallest size of an $s$-fold blocking set w.r.t. $(v-1)$-spaces in $\PG(r-1,q)$ is denoted by $b_q(r,h,r-v;s)$. Whenever the parameters are clear from the context, we just speak of generalized blocking sets.
 
A first attempt to generalize the Solomon--Stiffler construction is given by the following lemma.
\begin{lemma}
  \label{lemma_solomon_stiffler_gen}
  Let $h,f,r\in\N$ with $h+f\le r$, $\varepsilon_i\in\N$ for $h+f\le i\le r-1$ and $\sigma\in\N$ sufficiently large, e.g.\ $\sigma\ge\sum_{i=h+f}^{r-1}\varepsilon_i$. Then, we have $n_q(r,h,f;s)\ge n$, where $n:=\sigma\cdot\qbin{r}{h}{q}-\sum_{i=h+f}^{r-1}\varepsilon_i\cdot\qbin{i}{h}{q}$ and $s:=\sigma\cdot \qbin{r-f}{h}{q}-\sum_{i=h+f}^{r-1}\varepsilon_i\cdot\qbin{i-f}{h}{q}$.
\end{lemma}
\begin{proof}
  Consider the following multiset $\cM$ of $(h-1)$-spaces in $\PG(r-1,q)$. Starting from $\sigma$ copies of every $(h-1)$-space in $\PG(r-1,q)$ we remove the $(h-1)$-spaces contained in $\varepsilon_i$ $(i-1)$-spaces for all $h+f\le i\le r-1$, so that $|\cM|=n$. Since the set of all $(h-1)$-spaces contained in an arbitrary $(i-1)$-space is an $\qbin{i-f}{h}{q}$-fold blocking set w.r.t.\ $(r-f-1)$-spaces in $\PG(r-1,q)$, every codimension $f$ space in $\PG(r-1,q)$ contains at
  most $s$ elements from $\cM$.
\end{proof}
Choosing $\sigma=\varepsilon_4=1$ we e.g.\ obtain $n_2(5,2,2;6)\ge 120$. 
Similarly, $\sigma=\varepsilon_4=2$ yields $n_2(5,2,2;12)\ge 240$.

The essential idea in the proof of Lemma~\ref{lemma_solomon_stiffler_gen} is the blocking property of subspaces, so that we state the following alternative.
\begin{lemma}
  \label{lemma_remove_blocking_sets}
  Assume $h,f,r\in\N$ with $h+f\le r$, $l\in\N$, $\varepsilon_i\in\N$ for $1\le i\le l$. Let $\mathcal{B}_i$ be an $s_i$-fold blocking set of $(h-1)$-spaces with respect to $(r-f-1)$-spaces for every $1\le i\le l$. Moreover, let $|\mathcal{B}_i| = n_i$, and let $\sigma\in\N$ be sufficiently large. Then, we have $n_q(r,h,f;s)\ge n$, where $n:=\sigma\cdot\qbin{r}{h}{q}-\sum_{i=1}^{l}\varepsilon_i\cdot n_i$ and $s:=\sigma\cdot \qbin{r-f}{h}{q}-\sum_{i=1}^{l}\varepsilon_i\cdot s_i$.
\end{lemma}
\begin{proof}
  Form a multiset $\cM$ of $(h-1)$-spaces in $\PG(r-1,q)$, starting from $\sigma$ copies of every $(h-1)$-space in $\PG(r-1,q)$, and removing the $(h-1)$-spaces contained in $\epsilon_i$ copies of $\mathcal{B}_i$.
  The value for $s$ follows from the definition of $s_i$-fold blocking set: every codimension $f$ subspace contains at least $s_i$ elements from $\mathcal{B}_i$, so removing $\epsilon_i$ copies reduces the count by $\epsilon_i \cdot s_i$.
\end{proof}

In Section~\ref{sec_blocking_sets} we will consider blocking sets that have a smaller cardinality than the one consisting of all $(h-1)$-spaces in a fixed subspace.

A simple but very effective variant of Lemma~\ref{lemma_remove_blocking_sets} is given by the removal of a single blocking set.
\begin{lemma}
  \label{lemma_remove_one_blocking_set}
  If $\cB$ is a $s$-fold blocking set of $(h-1)$-spaces in $\PG(r-1,q)$ with respect to subspaces of codimension $f$ that has maximum multiplicity at most $m$, then we have 
  \begin{equation}
    n_q\!\left(r,h,f;m\cdot\qbin{r-f}{h}{q}-s\right)\ge m\cdot \qbin{r}{h}{q}-\left|\cB\right|. 
  \end{equation}
\end{lemma}

\section{Blocking sets}
\label{sec_blocking_sets}

In Section \ref{sec_intro} we have introduced the notion $b_q(r,h,f;s)$ for the minimum size of an $s$-fold blocking set of $(h-1)$-spaces with respect to subspaces of codimension $f$. Recall that this is the minimum number of $(h-1)$-spaces in $\PG(r-1,q)$ such that each subspace of codimension $f$ contains at least $s$ members. 
If we additionally assume that the maximum multiplicity of an $(h-1)$-space is $m$, then we use the notation $b_q(r,h,f;s,m)$ for the minimum possible cardinality. So, we obviously have $b_q(r,h,f;s,m)\ge b_q(r,h,f;s,m')$ and $b_q(r,h,f;s,m)\ge b_q(r,h,f;s)$ for all $m,m'\in\N$ with $m\le m'$. 
We note that allowing subspaces of projective dimension smaller than $h-1$ would not decrease those numbers. A straightforward counting argument gives a first lower bound. 

\begin{lemma}
\label{lemma_one_weight_bound_blocking_set}
We have
\begin{equation}
  b_q(r,h,f;s)\ge \frac{\qbin{r}{f}{q}\cdot s}{\qbin{r-h}{f}{q}}
  =
  \prod_{i=0}^{f-1} \frac{[r-i]_q}{[r-h-i]_q}\cdot s.
\end{equation}
\end{lemma}
\begin{proof}
   Each $(h-1)$-dimensional element is contained in
  $\qbin{r-h}{f}{q}$ subspaces of codimension $f$ and there are $\qbin{r}{f}{q}$ subspaces of codimension $f$ in total.
  This implies $ b_q(r,h,f;s)\ge \tfrac{\qbin{r}{f}{q}\cdot s}{\qbin{r-h}{f}{q}}$.
\end{proof}

A well-known construction is to use all $(h-1)$-spaces contained in a fixed $(h+f-1)$-space. Furthermore, note that the union of two blocking sets again gives a blocking set. Using this, we get the following proposition and lemma.

\begin{proposition}
  \label{prop_subspace_blocking}
  For $r,h,f\in\N$ with $h+f\le r$
  we have $b_q(r,h,f;1)\le \qbin{h+f}{h}{q}$.
\end{proposition}

\begin{lemma}
  \label{lemma_blocking_union}
  $
    b_q\!\left(r,h,f;s_1+s_2\right)\le
    b_q\!\left(r,h,f;s_1\right)+b_q\!\left(r,h,f;s_2\right)
  $.
\end{lemma}  

From Proposition~\ref{prop_subspace_blocking} we can e.g.\ conclude $b_q(5,2,2;1)\le \qbin{4}{2}{q}=q^4+q^3+2q^2+q+1$. Next we describe an improved construction from \cite{eisfeld1997blocking}. 
For a given integer $l\ge 3$ consider $\PG(2l-2,q)$ and an arbitrary point $P$.  With this, let  $\cS$ be a set of $\frac{q^{2l-2}-1}{q^2-1}$ planes through $P$ that  form a geometric line spread in the quotient geometry through $P$. Fix an $l$-space $U$ through $P$.  Let
$\cB$ consist of the $\frac{q^l-1}{q-1}$ lines in $U$ through $P$ together with the $q^2\frac{q^{2l-2}-1}{q^2-1}$ lines that lie in a plane of $\cS$ but do not contain $P$. Then every $(l-1)$-space in $\PG(2l-2,q)$ contains at least one line from $\cB$.

This construction is shown to be optimal in the theorem below. 

\begin{theorem}(\cite[Theorem 1.2]{eisfeld1997blocking})
  \label{thm_eisfeld_characterization}
  For each integer $l\ge 3$ we have 
  $b_q(2l-1,2,l-1;1)\ge \frac{q^{2l}-q^2}{q^2-1}+\frac{q^l-1}{q-1}$ and the above example is the only one in which equality holds.
\end{theorem}
We remark that in \cite{eisfeld1997blocking} sets of lines were considered. However, the statement remains obviously true for multisets of lines. For $l=3$ we obtain $b_q(5,2,2;1)= q^4+2q^2+q+1$, i.e.\ the subspace construction is improved by $q^3$ lines.

\subsection{The minimum number of lines in $\PG(4,q)$ such that every plane contains at least $s$ elements}
\label{subsec_bs_5_2_3}

In this subsection we want to focus on the values $b_q(5,2,2;s)$. From Theorem~\ref{thm_eisfeld_characterization}, Proposition~\ref{prop_subspace_blocking}, and Lemma~\ref{lemma_blocking_union} we directly conclude:
\begin{lemma}\label{lemblocksum}
  For each $0\le s'\le q^2+q$ and $t\ge 0$ with $s=t[3]_q+s'$ we have $b_q(5,2,2;s'+t[3]_q)\le \left(q^4+2q^2+q+1\right)\cdot s'+\qbin{5}{2}{q}\cdot t$. 
\end{lemma}
In the following we present some constructions that are better for specific choices for $s$ and we look into lower bounds improving upon Lemma~\ref{lemma_one_weight_bound_blocking_set}.

\begin{lemma}
  \label{lemma_q_fold_blocking_set_lines_vs_planes}
  In $\PG(4,q)$ there exists a $q$-fold blocking set w.r.t.\ planes consisting of $q^2(q+1)+ q^3 \left (q^2+1\right)$ (pairwise different) lines.
\end{lemma}
\begin{proof}
  For a point $P$, let $\cP_1,\dots,\cP_q$ be $q$ disjoint sets of $q^2+1$ planes, all containing the point $P$ and each forming a line spread in the factor geometry through $P$. Hence in the factor geometry, these line spreads are contained in a parallelism. With this, consider the set $\cL_1$ of lines consisting of all $q^2\cdot q\cdot \left(q^2+1\right)$ lines contained in one of the planes of the $\cP_i$ that do not contain $P$. Let $\cL_2$ be the $q^2(q+1)$ lines that do contain $P$ but are disjoint to  a fixed line $L$, which is contained in one of the planes in one of the $\cP_i$. Then we denote $\cL$ by $\cL_1\cup \cL_2$.
  Let $\pi$ be a plane with $P\in \pi$. If $\pi=\langle P, L \rangle$, then $\pi$ contains all $q^2$ lines of $\cL$ not through $P$. If $P\in \pi$, but $\pi\neq \langle P, L\rangle$, then $\pi$ intersects $\langle P, L\rangle$ in a point or in one line and hence $\pi$ contains at least $q$ lines from $\cL$ that contain point $P$.

  Let $\pi$ be a plane not containing $P$. The image of $\pi$ in the factor geometry through $P$ contains one line in each of the lines spreads, so that $\pi$ contains $q$ lines from $\cL$.
\end{proof}

For $q=2$ the corresponding blocking set has size $52$, which is indeed the minimum size for a double blocking set as verified by a small ILP computation, see Lemma~\ref{lemma_special_ilp_1}.
{In Lemma~\ref{lemma_bs_lb_1} we give a lower bound, which shows that $|B|\geq q^5+q^3+q^2+q$, so that Lemma~\ref{lemma_q_fold_blocking_set_lines_vs_planes} can be improved by at most $q^3-q$.} 
\begin{lemma}
  \label{lemma_qp1_fold_bs}
  In $\PG(4,q)$ there exists a $(q+1)$-fold blocking set w.r.t.\ planes consisting of $\qbin{4}{1}{q}+(q+1)q^2(q^2+1)$ (pairwise different) lines.
\end{lemma}
\begin{proof}
  Let $P$ be an arbitrary point in $\PG(4,q)$ and $\cP_1,\dots,\cP_{q+1}$ be sets of $q^2+1$ planes containing $P$, with the extra property that each forms a partial line parallelism in the factor geometry through $P$.
  With this, consider the set of lines $\cL$ consisting of all $q^2(q+1)\left(q^2+1\right)$ lines contained in one of the planes of the $\cP_i$ that do not contain $P$ and the $\qbin{4}{1}{q}$ lines that contain $P$. Denote the corresponding set of $\qbin{4}{1}{q}+(q+1)q^2(q^2+1)$ lines by $\cL$.

  Let $\pi$ be an arbitrary plane that contains $P$. The $q+1$ lines in $\pi$ that contain point $P$ are all contained in $\cL$.

  Let $\pi'$ be a plane not containing $P$. The image of $\pi'$ in the factor geometry through $P$ contains one line in each of the lines spreads, so that $\pi'$ contains $q+1$ lines from $\cL$.
\end{proof}

For $q=2$ this gives a $3$-fold blocking set of cardinality $75$.
By a sequence of ILP computations, see Lemma~\ref{lemma_special_ilp_2}, we can verify that this is the minimum possible cardinality (even for multisets of lines).

\begin{lemma}
  \label{lemma_q_2_fold_blocking_set}
  In $\PG(4,q)$ there exist a $q^2$-fold blocking set w.r.t.\ planes consisting of $q^6+q^4+q^3+q^2$ lines (maximum line multiplicity $q^2$).
\end{lemma}
\begin{proof}
  For an arbitrary plane $E$ let the multiset of lines $\cL$ consist of $q^2$ copies of each line in $E$ and a single copy of each of the $q^6$ lines outside of $E$, so that $|\cL|=q^6+q^4+q^3+q^2$.

  Now let $\pi$ be an arbitrary plane. If $\pi$ intersects $E$ in a line, then this line is contained $q^2$ times in $\cL$. For $\pi=E$ we have $q^2+q+1$ lines in $\cL$, each with multiplicity $q^2$. If $\pi$ intersects $E$ in a point, then $\pi$ contains $q^2$ lines disjoint to $E$, which are all contained in $\cL$.
\end{proof}
For $q=2$ this gives a $4$-fold blocking set of cardinality $92$, whose minimality can be concluded from Lemma~\ref{lemma_bs_lb_1}, see Theorem~\ref{thm_bs_5_2_2_q_2}. However, there are lines that are taken four times. In Appendix~\ref{sec_bs_lists} we list examples showing $b_2(5,2,2;4,1)\le 102$ and $b_2(5,2,2;4,2)\le 98$. The best known and indeed optimal construction for $b_2(5,2,2;5)$ is given by $b_2(5,2,2;5)\le b_2(5,2,2;1)+b_2(5,2,2;4)=27+92=119$, see Lemma~\ref{lemma_blocking_union}. Again, this construction comes with a large maximum line multiplicity. Moreover, using ILP, we found examples to prove $b_2(5,2,2;5,1)\le 123$, $b_2(5,2,2;5,3)\le 121$, and $b_2(5,2,2;5,4)\le 120$, see Appendix~\ref{sec_bs_lists}. 

\begin{lemma}
  \label{lemma_q_2_p_q_fold_blocking_set}
  In $\PG(4,q)$ there exist a $(q^2+q)$-fold blocking set w.r.t.\ planes consisting of $q^6+q^5+q^4+2q^3+2q^2+q=[6]_q+q^3+q^2-1$ lines (with maximum line multiplicity $q^2+q$).
\end{lemma}
\begin{proof}
  Let $L$ be an arbitrary line and $S\supseteq L$ be an arbitrary solid. With this, let $\cL_1$ be the set of all $(q+1)^2q$ lines that intersect $L$ in a point and are contained in $S$. Moreover, let $\cL_2$ be the set of lines that intersect $S$ in a point and are disjoint to $L$. As blocking set $\cB$ we choose $q^2+q$ times the line $L$, $q$ times the elements of $\cL_1$, and once the elements of $\cL_2$, so that $|B|=\left(q^2+q\right)+(q+1)^2q^2+(q+1)q^5=q^6+q^5+q^4+2q^3+2q^2+q$.

  Now we check the possible cases for a plane $\pi$. If $L\subseteq \pi$, then $L$ is contained $q^2+q$ times in $\cB$. If $\pi$ is disjoint to $L$, then $q^2+q$ elements of $\cL_2$ are contained in $\pi$. If $\pi$ intersects $L$ in a point and is contained in $S$, then $\pi$ contains $q+1$ elements from $\cL_1$. If $\pi$ intersects $L$ in a point and is not contained in $S$, then $\pi$ contains one element from $\cL_1$ and $q^2$ elements from $\cL_2$.
\end{proof}

For $q=2$ Lemma~\ref{lemma_q_2_p_q_fold_blocking_set} gives a $6$-fold blocking set of cardinality $138$, whose optimality is implied by Lemma~\ref{lemma_bs_lb_1}, see Theorem~\ref{thm_bs_5_2_2_q_2_pq}. However, there exists a line that is taken six times, so that we give examples showing $b_2(5,2,2;6,1)\le 146$, $b_2(5,2,2;6,2)\le 142$, $b_2(5,2,2;6,3)\le 142$, and $b_2(5,2,2;6,5)\le 141$ in Appendix~\ref{sec_bs_lists}.

\begin{table}[htp]
  \begin{center}
    \begin{tabular}{rrrrr}
      \hline
      $s$ & $b_2(5,2,2;s)$ & construction & lower bound & $b_2(5,2,2;s,1)\le$\\
      \hline
      1  & 27 & Theorem~\ref{thm_eisfeld_characterization} & Theorem~\ref{thm_eisfeld_characterization} & 27\\
      2  & 52 & Lemma~\ref{lemma_q_fold_blocking_set_lines_vs_planes}  & Lemma~\ref{lemma_special_ilp_1} & 52 \\
      3  & 75 & Lemma~\ref{lemma_qp1_fold_bs} & Lemma~\ref{lemma_special_ilp_2} & 75 \\
      4  & 92 & Lemma~\ref{lemma_q_2_fold_blocking_set} & Lemma~\ref{lemma_bs_lb_1} & 98 \\
      5  & 119 & Lemma~\ref{lemma_blocking_union} & Lemma~\ref{lemma_special_ilp_3} & 123 \\
      6  & 138 & Lemma~\ref{lemma_q_2_p_q_fold_blocking_set} & Lemma~\ref{lemma_bs_lb_1} & 146 \\
      7  & 155 & Proposition~\ref{prop_subspace_blocking} & Lemma~\ref{lemma_one_weight_bound_blocking_set} & 155\\
      \hline
    \end{tabular}
    \caption{Exact values for $b_2(5,2,2;s)$ and upper bounds for $b_2(5,2,2;s,1)$.}
    \label{table_b_2_5_2_2}
  \end{center}
\end{table}

In Table~\ref{table_b_2_5_2_2} we have summarized the upper bounds for $b_2(5,2,2;s)$ based on the constructions described so far, where $1\le s\le 7$. For future reference we have added the currently best known upper bound for $b_2(5,2,2;s,1)$ in the last column. Either the mentioned construction in the unrestricted cases automatically satisfies a maximum line multiplicity of one or the example was found by ILP computations. In the remaining part of this subsection we will present matching lower bounds. The lower bounds for $s\in\{2,3,5\}$ were obtained by tailored ILP computations, see Lemma~\ref{lemma_special_ilp_1}, Lemma~\ref{lemma_special_ilp_2}, and Lemma~\ref{lemma_special_ilp_3}, using the \texttt{ILOG CPLEX} solver.

\begin{lemma}
  \label{lemma_special_ilp_1}
  $b_2(5,2,2;2)\ge 52$.
\end{lemma}
\begin{proof}
  Direct ILP computation.
\end{proof}

\begin{lemma}
  \label{lemma_special_ilp_2}
  $b_2(5,2,2;3)\ge 75$.
\end{lemma}
\begin{proof}
  We utilize several ILP computations. If the maximum line multiplicity is $3$, then the minimum possible cardinality is $75$. 
  If the maximum line multiplicity is $2$ and there are two lines $L$, $L'$ with multiplicity $2$, then the minimum possible cardinality is $75$, independent of $\dim(L\cap L')$. 
  If the maximum line multiplicity of $2$ is attained at a unique line, then the minimum possible cardinality is at least $75$. 
  If the maximum line multiplicity is one and there exists a plane with three contained lines through a point, then the minimum possible cardinality is $75$. 
  If the maximum line multiplicity is one then there has to be a configuration as described before see Lemma \ref{lemma_qp1_fold_bs}.
\end{proof}

\begin{lemma}
  \label{lemma_special_ilp_4}
  The unique example attaining $b_2(5,2,2;4)= 92$ is given by the construction in the proof of Lemma~\ref{lemma_q_2_fold_blocking_set}.
\end{lemma}
\begin{proof}
  We utilize several ILP computations. If the maximum line multiplicity is at most three, then the cardinality is larger than $92$. If there is a unique line with multiplicity $4$, then the cardinality is larger than $94$. If there are two disjoint lines with multiplicity four, then the minimum possible cardinality is $95$. So, we prescribe two intersecting lines $L$, $L'$ with multiplicity four and minimize the chosen number of lines in $E:=\left\langle L,L'\right\rangle$ given a cardinality of $92$. It turns out that all seven lines in $E$ need to have multiplicity $4$ each. Prescribing such a configuration and cardinality $92$ results in a unique ILP solution. 
\end{proof}

\begin{lemma}
  \label{lemma_special_ilp_3}
  $b_2(5,2,2;5)\ge 119$.
\end{lemma}
\begin{proof}
  We utilize several ILP computations. If there is a line $L$ with multiplicity at least $6$, then the minimum possible cardinality is $120$. For maximum line multiplicity five the minimum possible cardinality is $119$. 
  For maximum line multiplicity at most four we considered a pair of lines $L$, $L'$ intersecting in a point, where $L$ attains the maximum multiplicity and $L'$ has the largest possible multiplicity of all lines intersecting $L$. For each choice of these two multiplicities we have checked by an ILP computation that cardinality $118$ is infeasible. 
\end{proof}

\begin{lemma}
  \label{lemma_bs_lb_1}
  Let $\cB$ be an $s$-fold blocking set of lines w.r.t.\ planes in $\PG(4,q)$. Then, we have
  $$
    \left|\cB\right| \ge \left(q^4+q^2+q+1\right) \cdot s -q\left(q+1\right)\cdot \cB(L)
  $$
  for each line $L$, where $\cB(L)$ denotes its multiplicity in $\cB$.
\end{lemma}
\begin{proof}
Fix a line $L$ and let $\cP_1$ be the set of planes that contain $L$ and $\cP_2$ the set of planes that are disjoint to $L$, so that $|\cP_1|=q^2+q+1$ and $|\cP_2|=q^6$. Consider the multiset $\cP:=q^2\cdot\cP_1+\cP_2$ of $q^6+q^4+q^3+q^2$ planes. Note that $L$ is contained in all elements of $\cP_1$ and therefore in $q^2\cdot\left(q^2+q+1\right)$ elements of $\cP$. Furthermore, any line $L'$ that is disjoint to $L$ is contained in $q^2$ elements of $\cP_2$ and $\cP$. And all other lines (i.e.\ those that intersect $L$ in a point) are contained in a unique element from $\cP_1$ and therefore in $q^2$ elements from $\cP$. Consider an $s$-fold blocking set $\cB$ of lines with respect to planes.
  We double count the set $S=\{(l, \pi)\,:\, l \in \cB, l\subset  \pi, \pi\in \cP\}$; which gives that
  \begin{align*} 
     \cB (L)  q^2(q^2+q+1)+\sum_{l' \neq L, l'\cap L\neq \emptyset} \cB(l') q^2 + \sum_{ l'\cap L= \emptyset} \cB(l') q^2 \geq (q^6+q^4+q^3+q^2) s,
     \end{align*}
    which is equivalent to 
    \begin{align*}
      q^2|\cB| + \cB (L)  q^2(q^2+q) \geq (q^6+q^4+q^3+q^2) s,
  \end{align*} and hence, proves the lemma.
\end{proof}

We note that Lemma~\ref{lemma_line_mult_bound_1} is a complementary bound.

\begin{lemma}
  \label{lemma_workhorse_finite_bs}
  If $b_q(5,2,2;s)\le \left(q^4+q^2+q+1\right)\cdot s$, then we have $$b_q(5,2,2;s+t[3]_q)=b_q(5,2,2;s)+t\qbin{5}{2}{q},$$ for all $t\in\N$. 
\end{lemma}
\begin{proof}
  Since $b_q(5,2,2;[3]_q)=\qbin{5}{2}{q}$, Lemma~\ref{lemma_blocking_union} yields $b_q(5,2,2;s+t[3]_q)\le b_q(5,2,2;s)+t\qbin{5}{2}{q}$ for all $t\in\N$. Now assume that $\cB$ is a $(s+t[3]_q)$-fold blocking set of $n$ lines w.r.t.\ planes where $n<b_q(5,2,2;s)+t\qbin{5}{2}{q}$. W.l.o.g.\ we assume that $t$ is minimal with this property, which implies the existence of a line $L$ with multiplicity $\cB(L)=0$. Lemma~\ref{lemma_bs_lb_1} gives 
  \begin{eqnarray*}
    \left|\cB\right| &\ge& \left(q^4+q^2+q+1\right)\cdot (s+t[3]_q)\\ 
    &=& \left(q^4+q^2+q+1\right)\cdot s +t\qbin{5}{2}{q}+tq(q+1), 
  \end{eqnarray*}
  which is a contradiction.
\end{proof}

\begin{theorem}
  \label{thm_bs_5_2_2_q_2}
  For all $t\in\N$ we have 
  $$b_q\!\left(5,2,2;t[3]_q+q^2\right)=t\qbin{5}{2}{q}+q^2\cdot\left(q^4+q^2+q+1\right).$$ If equality is attained, then every line has multiplicity at least $t$.
\end{theorem}
\begin{proof}
  Lemma~\ref{lemma_q_2_fold_blocking_set} gives a matching construction for $t=0$, so that Proposition~\ref{prop_subspace_blocking} and Lemma~\ref{lemma_blocking_union} imply the corresponding upper bound for all $t\in\N$. Lemma~\ref{lemma_bs_lb_1} gives a matching lower bound for $t=0$, so that the statement follows from Lemma~\ref{lemma_workhorse_finite_bs}.
\end{proof}
Since the construction in the proof of Lemma~\ref{lemma_q_2_fold_blocking_set} gives the unique $4$-fold blocking set of $92$ lines in $\PG(4,2)$ w.r.t.~planes, see Lemma~\ref{lemma_special_ilp_4}, there is a unique example attaining $b_2(5,2,2;4+7t)$ for all $t\in \N$.

\begin{theorem}
  \label{thm_bs_5_2_2_q_2_pq}
  For all $t\in\N$ we have 
  $$b_q\!\left(5,2,2;t[3]_q+q^2+q\right)=t\qbin{5}{2}{q}+\left(q^2+q\right)\left(q^4+q^2+q+1\right).$$ If equality is attained, then every line has multiplicity at least $t$.
\end{theorem}
\begin{proof}
  Lemma~\ref{lemma_q_2_p_q_fold_blocking_set} gives a matching construction for $t=0$, so that Proposition~\ref{prop_subspace_blocking} and Lemma~\ref{lemma_blocking_union} imply the corresponding upper bound for all $t\in\N$. Lemma~\ref{lemma_bs_lb_1} gives a matching lower bound for $t=0$, so that the statement follows from Lemma~\ref{lemma_workhorse_finite_bs}.
\end{proof}

\begin{theorem}
  \label{thm_b_2_5_2_2}
  For each $t\in \N$ we have 
  \begin{itemize}
    \item $b_2(5,2,2;1+7t)=27+155t$,
    \item $b_2(5,2,2;2+7t)=52+155t$,
    \item $b_2(5,2,2;3+7t)=75+155t$,
    \item $b_2(5,2,2;4+7t)=92+155t$,
    \item $b_2(5,2,2;5+7t)=119+155t$,
    \item $b_2(5,2,2;6+7t)=138+155t$,
    \item $b_2(5,2,2;7+7t)=155(t+1)$.
  \end{itemize}
\end{theorem}
\begin{proof}
  For the constructions and upper bounds for $b_2(5,2,2;s)$ for $1\le s\le 7$ we refer to Table~\ref{table_b_2_5_2_2}. 
  Since $n_2(5,2,2;7)=155$, Lemma~\ref{lemma_blocking_union} extends the upper bounds to all $t\in\N$. Therefore, it remains to give the lower bounds.
  Lemma~\ref{lemma_bs_lb_1} with $\cB(L)=0$ shows that these upper bounds for $b_2(5,2,2;s)$ are tight for $s\in\{4,6\}$, so that we can apply Lemma~\ref{lemma_workhorse_finite_bs}. Applying Lemma~\ref{lemma_bs_lb_1} with $s=8$ and $\cB(L)=0$ would give a lower bound of $184>182$, so that we may suppose that $\cB(L)\ge 1$ for each $L$.  Hence, it suffices to determine $b_2(5,2,2;1)$, which is done in \cite{eisfeld1997blocking}. Applying Lemma~\ref{lemma_bs_lb_1} with $\cB(L)=0$ gives a lower bound that matches the size of the stated constructions for $s\in\{9,10\}$, see Table~\ref{table_b_2_5_2_2} and Lemma \ref{lemma_blocking_union}, and is strictly larger for $s=12$. So, it suffices to determine $b_2(5,2,2;s)$ for $s\in\{2,3,5\}$, see Lemma~\ref{lemma_special_ilp_1}, Lemma~\ref{lemma_special_ilp_2}, and Lemma~\ref{lemma_special_ilp_3} for the corresponding lower bounds.
\end{proof}

In Table~\ref{table_b_3_5_2_2} we fix $q=3$ and summarize our knowledge on $b_3(5,2,2;s)$ for $1\le s\le 13$.

\begin{table}[htp]
  \begin{center}
    \begin{tabular}{rrrr}
      \hline
      $s$ & $b_3(5,2,2;s)$ & construction & lower bound\\
      \hline
      1  & 103 & Theorem~\ref{thm_eisfeld_characterization} & Theorem~\ref{thm_eisfeld_characterization} \\
      2  & 188--206 & Lemma~\ref{lemma_blocking_union}  & Lemma~\ref{lemma_bs_lb_1} \\
      3  & 282--306 & Lemma~\ref{lemma_q_fold_blocking_set_lines_vs_planes} & Lemma~\ref{lemma_bs_lb_1} \\
      4  & 376--400 & Lemma~\ref{lemma_qp1_fold_bs} & Lemma~\ref{lemma_bs_lb_1} \\
      5  & 470--502 & ILP & Lemma~\ref{lemma_bs_lb_1} \\
      6  & 564--600 & ILP & Lemma~\ref{lemma_bs_lb_1} \\
      7  & 658--690 & ILP & Lemma~\ref{lemma_bs_lb_1} \\
      8  & 752--784 & ILP & Lemma~\ref{lemma_bs_lb_1} \\
      9  & 846 & Lemma~\ref{lemma_q_2_fold_blocking_set}  & Lemma~\ref{lemma_bs_lb_1} \\
      10 & 940--949 & Lemma~\ref{lemma_blocking_union} & Lemma~\ref{lemma_bs_lb_1} \\
      11 & 1034--1050 & ILP & Lemma~\ref{lemma_bs_lb_1} \\
      12 & 1128 & Lemma~\ref{lemma_q_2_p_q_fold_blocking_set} & Lemma~\ref{lemma_bs_lb_1} \\
      13 & 1210 & Proposition~\ref{prop_subspace_blocking} & Lemma~\ref{lemma_one_weight_bound_blocking_set} \\
      \hline
    \end{tabular}
    \caption{Bounds  for $b_3(5,2,2;s)$.}
    \label{table_b_3_5_2_2}
  \end{center}
\end{table}

\subsection{Generalizations to other parameters}
\label{subsec_bs_generalizations}

The constructions from Lemma~\ref{lemma_q_2_fold_blocking_set} and Lemma~\ref{lemma_q_2_p_q_fold_blocking_set} can be described from a more general point of view. 
In $\PG(r-1,q)$ let $S_1,\dots,S_{r-1}$ be a chain of subspaces with $\dim\!\left(S_i\right)=i-1$ for $1\le i\le r-1$. 
The set of $(h-1)$-spaces is partitioned into classes $\cH_1,\dots,\cH_u$ according to the intersection dimensions with those $S_i$. The set of subspaces of codimension $f$ is partitioned into classes $\cF_1,\dots,\cF_v$ according to the intersection dimensions with the $S_i$. 
By $\beta_{i,j}$ we denote the number of elements from $\cH_i$ that are contained in an arbitrary element $\pi\in\cF_j$. As blocking set we choose $\cB=\sum_{i=1}^u \alpha_i\cH_i$, where $\alpha_i\in\N$ for $1\le i\le u$.  Given this framework, we obtain a simple optimization problem: choose $\alpha_i\in \N$ minimizing $\sum_{i=1}^u \alpha_i\cdot \left|\cH_i\right|$ such that $\sum_{i}\alpha_i\beta_{i,j}\ge s$ for all $1\le j\le v$. Of course also the easy construction from Proposition~\ref{prop_subspace_blocking} can be described in this way.

\begin{table}[htp]
  \begin{center}
  \begin{tabular}{rcrrcr}
  \hline
  $i$ & $\cH_i$ & $\#$ & $j$ & $\cF_j$ & $\#$ \\
  \hline
  1 & $(0,1,1,1)$ & $1$ & 1 & $(0,1,2,2)$ & $1$ \\
  2 & $(0,0,1,1)$ & $q$ & 2 & $(0,1,1,2)$ & $q$ \\
  3 & $(-1,0,1,1)$ & $q^2$ & 3 & $(0,0,1,2)$ & $q^2$ \\
  4 &$(0,0,0,1)$ & $q^3$ & 4 & $(-1,0,1,2)$ & $q^3$ \\
  5 & $(-1,0,0,1)$ & $q^3$ & 5 & $(0,1,1,1)$ & $q^2$ \\
  6 & $(-1,-1,0,1)$ & $q^4$ & 6 & $(0,0,1,1)$ & $q^3$ \\
  7 & $(0,0,0,0)$ & $q^3$ & 7 & $(-1,0,1,1,)$ & $q^4$ \\
  8 & $(-1,0,0,0)$ & $q^4$ & 8 &$(0,0,0,1)$ & $q^4$ \\
  9 & $(-1,-1,0,0)$ & $q^5$ & 9 & $(-1,0,0,1)$) & $q^5$ \\
  10 & $(-1,-1,-1,0)$ & $q^6$ & 10 & $(-1,-1,0,1)$ & $q^6$ \\
  \hline
  \end{tabular}
  \caption{Line and plane classes in $\PG(4,q)$ according to the intersection dimensions with a chamber, i.e.\ a maximal flag.}
  \label{table_classes}
  \end{center}
\end{table}

\begin{example}
  Let $K$ be a chamber in $\PG(4,q)$, which is a maximal flag $\{ \pi_0, \pi_1, \pi_2,\pi_3\}$, where $\pi_0\subset \pi_1\subset \pi_2\subset \pi_3$ and $\dim(\pi_i)=i$. For $\PG(4,q)$ and $(h,f)=(2,2)$ we obtain the ten line classes and ten plane classes listed in Table~\ref{table_classes}, according to their intersection dimensions $K$.
  Therefore, as an example, $\cH_3$ consists of all lines, contained in $\pi_2$ (and hence also in $\pi_3$) and meeting the line $\pi_1$ precisely in the point $\pi_0$. On the other hand, $\cF_8$ consists of all planes in $\pi_3$, that meet $\pi_2$ precisely in the point $P$. \\ The corresponding intersection numbers $\beta_{i,j}$ are given in Table~\ref{table_intersection}. Note that we have $\beta_{i,j}\in\!\left\{0,1,q,q^2\right\}$.
\end{example}

\begin{table}[htp]
  \begin{center}
  \begin{tabular}{ccccccccccc}
  \hline
  $j/i$ & 1 & 2 & 3 & 4 & 5 & 6 & 7 & 8 & 9 & 10\\ 
  \hline
  1&1&$q$&$q^2$&0&0&0&0&0&0&0\\ 
  2&1&0&0&$q$&$q^2$&0&0&0&0&0\\ 
  3&0&1&0&$q$&0&$q^2$&0&0&0&0\\ 4&0&0&1&0&$q$&$q^2$&0&0&0&0\\ 
  5&1&0&0&0&0&0&$q$&$q^2$&0&0\\ 
  6&0&1&0&0&0&0&$q$&0&$q^2$&0\\ 
  7&0&0&1&0&0&0&0&$q$&$q^2$&0\\ 
  8&0&0&0&1&0&0&$q$&0&0&$q^2$\\ 
  9&0&0&0&0&1&0&0&$q$&0&$q^2$\\ 
 10&0&0&0&0&0&1&0&0&$q$&$q^2$\\
  \hline
  \end{tabular}
  \caption{Intersection numbers $\beta_{i,j}$ in $\PG(4,q)$ w.r.t.\ lines and planes.}
  \label{table_intersection}
  \end{center}
\end{table}

\medskip

We can generalize Lemma~\ref{lemma_bs_lb_1} as follows.
\begin{lemma}
  \label{lemma_bs_lb_2}
  Let $h=2$, $f\ge h$, and $r>h+f$. Then, for any $s$-fold blocking set of lines in $\PG(r-1,q)$ w.r.t. to subspace of codimension $f$ we have
  $$
    \left|\cB\right| \ge \frac{1}{\beta_{2,1}}\cdot\left(\alpha_1+\frac{\beta_{2,1}-\beta_{3,1}}{\beta_{3,2}}\cdot \alpha_2\right) \cdot s -\frac{\beta_{1,1}-\beta_{2,1}}{\beta_{2,1}}\cdot \cB(L)
  $$
  for each line $L$, where $\cB(L)$ denotes its multiplicity in $\cB$, 
  \begin{eqnarray*}
    \alpha_1 &=& \qbin{r-2}{f}{q}, \qquad
    \alpha_2 = q^{2(r-f)}\cdot\qbin{r-2}{r-f}{q},\\
    \beta_{i,1} &=& \qbin{r-i-1}{f}{q} \qquad \text{for $i\in \{1,2,3\}$}\\
    \beta_{3,2} &=& q^{2(r-f-2)}\qbin{r-4}{f-2}{q}.
  \end{eqnarray*}
\end{lemma}
\begin{proof}
  Let $\cL_2$ be the set of lines in $\PG(r-1,q)$ that intersect $L$ in a point and $\cL_3$ be the set of lines that are disjoint to $L$. Set $\cL_1:=\{L\}$ and $\cL:=\cL_1\cup\cL_2\cup\cL_3$, i.e.\ the set of all lines in $\PG(r-1,q)$. By $\cP_1$ we denote the set of subspaces of codimension $f$ that contain $L$ and by $\cP_2$ we denote the set of subspaces of codimension $f$ that are disjoint to $L$. 
  
  Recall that the number of $j$-spaces disjoint from a fixed $m$-space in $\PG(n, q)$ equals $q^{(m+1)(j+1)}\qbin{n-m}{j+1}{q}$.
  With this we have
  \begin{equation}
    \alpha_1:=\left|\cP_1\right|=\qbin{r-2}{f}{q}
  \end{equation}
  and
  \begin{equation}
    \alpha_2:=\left|\cP_2\right|
    =q^{2(r-f)}\qbin{r-2}{r-f}{q}.
  \end{equation}
  {A line $l'\in \mathcal{L}_i$ is contained in an element $\pi\in \cP_1$ if the $i$-space $\langle l',L\rangle$ is contained in $\pi$. Hence, $\beta_{i,1} = \qbin{r-i-1}{f}{q}$. A line $l'\in \cL_1\cup \cL_2$ meets the line $L$ and hence, cannot be contained in an element of $\cP_2$; which implies $\beta_{1, 2} = \beta_{2, 2} = 0$. For a line $l'\in \cL_3$, we need to define the number $\beta_{3,2}$ of $(r-f-1)$-spaces in $\PG(r-1,q)$ through $l'$ and disjoint from $L$. This equals the number of $(r-f-3)$-spaces in $\PG(r-3,q)$, disjoint from a line, which equals $q^{2(r-f-2)}\qbin{r-4}{f-2}{q}$.}
  
  Let $\cB$ be an $s$-fold blocking set of lines in $\PG(r-1,q)$ w.r.t.\ a subspace of codimension $f$. Choose ${t\in\mathbb{R}_{\ge 0}}$ such that $\beta_{3,1}+t\cdot\beta_{3,2}=\beta_{2,1}$, i.e.\ 
  \begin{equation}
    t:=\frac{\beta_{2,1}-\beta_{3,1}}{\beta_{3,2}}.
  \end{equation}
  With this, we double count the set $S=\{(l,\alpha): l\in \cB, l\subset \alpha, \dim(\alpha)=r-f-1\}$.
  \begin{equation}
    \sum_{E\in\cP_1}\sum_{U\subseteq E\,:\,U\in\cL} \cB(U) \,+\, 
    t\cdot \sum_{E\in\cP_2}\sum_{U\subseteq E\,:\,U\in\cL} \cB(U) 
    \ge \left(\alpha_1+t\cdot\alpha_2\right)\cdot s.
  \end{equation}
  Note that $\cB(U)$ is counted $\beta_{2,1}$ times in the sum on the left hand side for all $U\in\cL_2\cup\cL_3$ while $\cB(L)$ is counted $\beta_{1,1}$ times.
\end{proof}

It seems tempting to generalize Theorem~\ref{thm_bs_5_2_2_q_2} (or Theorem~\ref{thm_bs_5_2_2_q_2_pq}). However, there are some issues that we cannot resolve. Motivated by Lemma~\ref{lemma_q_2_fold_blocking_set} we state the following generalized construction:
\begin{lemma}
\label{lemma_q_2_fold_blocking_set_gen}
  For $r\ge 4$ there exist a $q^2$-fold blocking set 
  in $\PG(r,q)$  w.r.t.\ planes consisting of $q^{2(r-1)}+q^2\cdot\qbin{r-1}{2}{q}$ lines (maximum line multiplicity $q^2$).
\end{lemma}
\begin{proof}
  For an arbitrary but fixed $(r-2)$-space $S$, let $\cL_1$ be the set of $\qbin{r-1}{2}{q}$ lines contained in $S$ and let $\cL_2$ the set of $q^{2(r-1)}$ lines disjoint to $S$. With this we set $\cB=q^2\cdot \cL_1+\cL_2$ and check that $\cB$ is indeed a $q^2$-fold blocking set w.r.t.\ planes.
\end{proof}
So, we especially have $b_q\!\left(6,2,3;q^2\right)\le q^8 + q^6+q^5+2q^4+q^3+q^2$. Applying Lemma~\ref{lemma_bs_lb_2} for these parameters with $\cB(L)=0$ gives $b_q\!\left(6,2,3;q^2\right)\ge q^8 + q^6+q^5+q^4+q^3+q^2$, i.e.\ there remains a gap of $q^4$. For $q=2$ those bounds give $380\le b_2(6,2,3;4)\le 396$. Solving our standard ILP model with the additional constraint $x_L=0$, i.e.\ $\cB(L)=0$, gives $b_2(6,2,3;4)=396$, so that Lemma~\ref{lemma_q_2_fold_blocking_set_gen} is optimal for $(r,q)=(5,2)$. We remark that the corresponding LP relaxation yields the lower bound     $b_2(6,2,3;4)\ge 380$ only, and hence, the bound in Lemma~\ref{lemma_bs_lb_2} is optimal for these parameters if we only rely on counting arguments and the extra information $\cB(L)=0$. Using $\cB(L)=0$ and $\cB(L')=0$ for two disjoint lines the corresponding LP relaxation gives $b_2(6,2,3;4)\ge 385.3333$. For $b_q\!\left(6,2,3;q^2+q\right)$ similar computations can be performed.

\section{The maximum number of lines in $\mathbf{\PG(4,q)}$ such that each plane contains at most $\mathbf{s}$ lines}
\label{sec_n_q_4_2_2_s}

Here we want to determine bounds for $n_q(5,2,2;s)$. Recall the correspondence, by duality, between $n_q(5,2,2;s)$ and the generalized blocking sets. More precisely, if we have a configuration with maximum line multiplicity at most $s$, then $n_q(5,2,2;s) = s\cdot \qbin{5}{2}{q}-b_q(5,2,2;s\cdot[3]_q-s,s)$.

\begin{theorem}
  \label{thm_n_q_5_2_2_1}
  We have $n_q(5,2,2;1)=q^3+1$.
\end{theorem}
\begin{proof}
  From Lemma~\ref{lemma_connection_to_subspace_codes} we conclude $n_q(5,2,2;1)=A_q(5,2;4)$. Here $A_q(5,2;4)$ is the maximum cardinality of a partial line spread in $\PG(4,q)$, which is well known, see e.g.\ \cite[Theorem 5]{nuastase2017maximum}.
\end{proof}

\begin{corollary}
    We have $b_q(5,2,2;q^2+q,1)=\qbin{5}{2}{q}-\left(q^3+1\right)$.
\end{corollary}
From Lemma \ref{thm_bs_5_2_2_q_2_pq} we know that $b_q(5,2,2;q^2+q) = (q^2+q)(q^4+q^2+q+1)$, which shows that here again, the last parameter $m$ plays an important role.

We can easily formulate the problem of the determination of $n_q(r,h,f;s)$ as an integer linear programming (ILP) problem. To this end let $\mathcal{H}$ denote the set of all $(h-1)$-spaces and $\mathcal{F}$ denote the set of all $(r-f-1)$-spaces in $\PG(r-1,q)$. As variables we choose $x_H\in\N$ for all $H\in\mathcal{H}$ to model the multiplicities of the chosen $(h-1)$-spaces. The condition that each subspace of codimension $f$ contains at most $s$ elements can be modeled as $\sum_{H\in\mathcal{H}\,:\, H\subseteq F} x_H\le s$ for all $F\in\mathcal{F}$. As target function we choose $\sum_{H\in\mathcal{H}} x_H$, i.e.\ the number of selected $(h-1)$-spaces. Typically this ILP can be solved directly for rather small values of $r$, $h$, $f$, and $q$ only. In order to obtain lower bounds we can e.g.\ prescribe some automorphisms. For upper bounds we can add tailored extra constraints or prescribe a few $(h-1)$-spaces to reduce the symmetry of the formulation.

From Theorem~\ref{thm_n_q_5_2_2_1} and Lemma~\ref{lemma_sum_gen_ham} we have $n_q(5,2,2;2)\ge 2\left(q^3+1\right)\in \Theta\!\left(q^3\right)$.
 From Lemma~\ref{lemma_one_weight_bound_gen_hamming_weight} we conclude $n_q(5,2,2;2)\le 2\cdot \frac{\left(q^4+q^3+q^2+q+1\right)\cdot\left(q^2+1\right)}{q^2+q+1}\in\Theta\!\left(q^4\right)$, so that the question for the right order of magnitude, in terms of $q$, arises. By ILP computations we found examples showing $n_2(5,2,2;2)\ge 32$, $n_3(5,2,2;2)\ge 97$, and $n_5(5,2,2;2)\ge 493$. 
We remark that improved constructions for $n_q(5,2,1;2)$ have been recently obtained in \cite{krotovkurz2025}.

\begin{lemma}
  \label{lemma_line_mult_bound_1}
  Let $\cL$ be a multiset of lines in $\PG(4,q)$ such that each
  plane contains at most $s$ lines. Then, we have
  \begin{equation}
    \left|\cL\right| \le \left(q^4+q^2+q+1\right)\cdot s-q(q+1)\cdot\cL(L),
  \end{equation}
  for each line $L$, where $\cL(L)$ denotes the multiplicity of $L$ in $\mathcal{L}$.
\end{lemma}
\begin{proof} 
Fix a line $L$ and let $\cP_1$ be the set of planes that contain $L$ and $\cP_2$ the set of planes that are disjoint to $L$. Then we have $|\cP_1|=q^2+q+1$ and $|\cP_2|=q^6$. Consider the multiset $\cP:=q^2\cdot\cP_1+\cP_2$ of $q^6+q^4+q^3+q^2$ planes. Note that $L$ is contained in all elements of $\cP_1$ and so in $q^2\cdot\left(q^2+q+1\right)$ elements of $\cP$. Any line $L'$ that is disjoint to $L$ is contained in $q^2$ elements of $\cP_2$ and $\cP$. Moreover, all other lines (i.e.\ those that intersect $L$ in a point) are contained in a unique element from $\cP_1$ and so $q^2$ elements from $\cP$. Consider a projective $(2,2)-(|\cL|,5,s)$ system $\cL$.

  We double count the set $S=\{(l, \pi)\,:\, l \in \cL, l\subset  \pi, \pi\in \cP\}$; which gives that
  \begin{align*} 
     \cL (L)  q^2(q^2+q+1)+\sum_{l' \neq L, l'\cap L\neq \emptyset} \cL(l') q^2 + \sum_{ l'\cap L= \emptyset} \cL(l') q^2 \leq (q^6+q^4+q^3+q^2) s,
     \end{align*}
    which is equivalent to 
    \begin{align*}
      q^2|\cL| + \cL (L)  q^2(q^2+q) \leq (q^6+q^4+q^3+q^2) s,
  \end{align*} and hence, proves the lemma.
\end{proof}

\begin{lemma}\label{lem_new}
    Let $B$ be the smallest $s$-fold blocking set of lines in $\PG(4,q)$ with respect to planes and maximum multiplicity $m$. Then 
    \begin{align*}
        155m-|B| = 155m-b_2(5,2,2;s,m) \leq n_2(5,2,2;7m-s) \leq 23(7m-s)-6m.
    \end{align*}
\end{lemma}
\begin{proof}
    Follows immediately from Lemmas \ref{lemma_remove_one_blocking_set} and  \ref{lemma_line_mult_bound_1}.
\end{proof}

\begin{lemma}
  We have $166\le n_2(5,2,2;8)\le 172$, $323\le n_2(5,2,2;15)\le 327$, and $478\le n_2(5,2,2;22)\le 482$.
\end{lemma}
\begin{proof}
    The upper bound follows immediately from Lemma \ref{lem_new} with $s=6$ and $m\in\{2,3,4\}$. For the lower bound, we also use the bounds $b_2(5,2,2;6,2)\le 144$,  $b_2(5,2,2;6,3)\le 142$ {and $b_2(5,2,2;6,4)\leq 142$}.
\end{proof}

\begin{lemma}
  We have $32\le n_2(5,2,2;2)\le 34$, 
  $187\le n_2(5,2,2;9)\le 195$, $344\le n_2(5,2,2;16)\le 350$, and $500\le n_2(5,2,2;23)\le 505$.
\end{lemma}
\begin{proof}
    The upper bound follows immediately from Lemma \ref{lem_new} for $s=5$ and $m\in\{1,2,3,4 \}$. For the lower bound, we also use the bounds $b_2(5,2,2;5,2)\le b_2(5,2,2;5,1)\le 123$,  $b_2(5,2,2;5,3)\le 121$ {and $b_2(5,2,2;5,4)\leq 120$}. In the case of $m=1$ we can find a better upper bound: let $\cP$ be a faithful $(2,2)-(n,5,2)_2$ projective system. If there exists a line $L$ in $\cP$ with multiplicity $\cP(L)$ at least $2$, then Lemma~\ref{lemma_line_mult_bound_1} implies $n\le 34$. For maximum line multiplicity one, we utilize an ILP computation to verify $n\le 34$.  
    $b_2(5,2,2;5,4)\le 120$.
\end{proof}

\begin{lemma}
  We have $53\le n_2(5,2,2;3)\le 59$, 
  $212\le n_2(5,2,2;10)\le 218$, $367\le n_2(5,2,2;17)\le 373$, and $n_2(5,2,2;24)=528$.
\end{lemma}
\begin{proof}
    The upper bound follows from Lemma \ref{lem_new} for $s=4$ and $m\in\{1,2,3,4 \}$. 
    For the lower bound, we also use the bounds $b_2(5,2,2;4,1)\le 102$, $b_2(5,2,2;4,2)\le 98$, and  $b_2(5,2,2;4,4)\leq 92$.  
    In the case of $m=1$ we can find a better upper bound: let $\cP$ be a faithful $(2,2)-(n,5,3)_2$ projective system. If there exists a line $L$ in $\cP$ with multiplicity $\cP(L)$ at least $2$, then Lemma~\ref{lemma_line_mult_bound_1} implies $n\le 57$. For maximum line multiplicity one, we utilize an ILP computation to verify $n\le 59$.

\end{proof}

\begin{proposition}
  For $t\in\N$ we have $n_2(5,2,2;7t+4)=155t+80$.    
\end{proposition}
\begin{proof}
From Lemma \ref{lem_new} with $s=3$ and $m=t+1$, we have that $155(t+1)-b_2(5,2,2; 3, 1)\leq 155(t+1)-b_2(5,2,2; 3, t+1)\leq n_2(5,2,2;7t+4)\leq 23(7t+4)-6\cL(l)$. Using $b_2(5,2,2;3,1)=b_2(5,2,2;3)=75$ we get the right lower bound. Now, let $\cP$ be a faithful $(2,2)-(n,5,7t+4)_2$ projective system. If there exists a line $L$ in $\cP$ with multiplicity $\cL(L)$ at least $t+2$, then we find the right lower bound $n\le 155t+80$. For maximum line multiplicity $t+1$ we conclude $n\le (t+1)\cdot 155-b_2(5,2,2;3,1)=155t+80$, which proves the statement.
\end{proof}

Actually, 
Lemma~\ref{lemma_remove_one_blocking_set} and the construction of a blocking set in Lemma~\ref{lemma_qp1_fold_bs} imply the following lemma.
\begin{lemma}
  \label{lemma_construction_5_2_2_q}
  We have $n_q(5,2,2;q^2)\ge q^4\cdot \left(q^2+1\right)$.
\end{lemma}
%

\begin{proposition}
  For $t\in\N$, we have $n_2(5,2,2;7t+5)=155t+ 103$.
\end{proposition}
\begin{proof}
From Lemma \ref{lem_new} with $s=2$ and $m=t+1$, we find the right lower bound using  $b_2(5,2,2;2,1)=52$. 
For the upper bound, let $\cP$ be a faithful $(2,2)-(n,5,7t+5)_2$ projective system. If there exists a line $L$ in $\cP$ with multiplicity $\cP(L)$ at least $t+2$, then Lemma~\ref{lemma_line_mult_bound_1} implies $n\le 155t+103$. For maximum line multiplicity $t+1$, we conclude $n\le (t+1)\cdot 155-b_2(5,2,2;2,1)=155t+103$.
\end{proof}


\begin{theorem}
  \label{thm_5_2_2_q_gen}
  For each $t\ge 0$, we have
  $$
    n_q\!\left(5,2,2;t\cdot[3]_q +q^2+q\right)=t\cdot\qbin{5}{2}{q}+q^6+q^5+q^4+2q^3.
  $$
\end{theorem}
\begin{proof}
  Consider the set of all lines in $\PG(4,q)$ with multiplicity $(t+1)$ and subtract those from a blocking set $\cB$ as in Theorem~\ref{thm_bs_5_2_2_q_2_pq}. Since
  the total number of lines is given by $\qbin{5}{2}{q}$, we have $n_q(5,2,2;t[3]_q+q^2+q)\ge (t+1)\cdot \qbin{5}{2}{q}-\left(q^4+2q^2+q+1\right)$.
  
  Now consider a multiset $\cL$ of lines in $\PG(4,q)$ such that each plane contains at most $t\cdot[3]_q+q^2+q$ lines and that $\left|\cL\right|>t\cdot\qbin{5}{2}{q}+q^6+q^5+q^4+2q^3$. If there exists a line $L$ with $\cL(L)\ge t+2$,
  then Lemma~\ref{lemma_line_mult_bound_1} yields
  \begin{eqnarray*}
    \left|\cL\right|&\le& \left(q^4+q^2+q+1\right)\cdot \left(t\cdot[3]_q+q^2+q\right)-(t+2)q(q+1)\\
    &=&t\cdot\qbin{5}{2}{q}+q^6+q^5+q^4+2q^3-q,
  \end{eqnarray*}
  which is a contradiction. Thus, the maximum line multiplicity $\cL(L)$ is at most $t$ and we denote the complementary multiset of lines by $\cB$. Since each plane contains at most $q^2+q=[3]_q-1$ elements from $\cS$, the elements of $\cB$ block every plane at least once. From  Theorem~\ref{thm_bs_5_2_2_q_2_pq} we conclude
  \begin{eqnarray*}
    \left|\cL\right|&=&(t+1)\cdot \qbin{5}{2}{q}-\left|\cB\right| \le
    (t+1)\cdot \qbin{5}{2}{q}-\left(q^4+2q^2+q+1\right)\\ 
    &=&t\qbin{5}{2}{q}+q^6+q^5+q^4+2q^3,
  \end{eqnarray*}
  which is a contradiction, and hence, proves the theorem.
\end{proof}

\medskip

We have the summarized our information on $n_2(5,2,2;s)$ 
in Table~\ref{table_bounds_n_2_5_2_2}. 

\begin{table}[htp]
  \begin{center}
    \begin{tabular}{rrrrrrrr}
      \hline
      $s$ & $n_2(5,2,2;s)$ & $s$ & $n_2(5,2,2;s)$ & $s$ & $n_2(5,2,2;s)$ & $s$ & $n_2(5,2,2;s)$\\
      \hline
      1&9     & 8&166--172&15&323--327&22&478--482\\
      2&32--34& 9&187--195&16&344--350&23&500--505\\
      3&53--59&10&212--218&17&367--373&24&528\\
      4&80    &11&235     &18&390&25&545\\
      5&103   &12&258     &19&413&26&568\\
      6&128   &13&283     &20&438&27&593\\
      7&155   &14&310     &21&465&28&620\\
      \hline
    \end{tabular}
    \caption{Bounds for $n_2(5,2,2;s)$.}
    \label{table_bounds_n_2_5_2_2}
  \end{center}
\end{table}

We can easily generalize Lemma~\ref{lemma_line_mult_bound_1} to $\PG(n,q)$. For an even more general version, formulated in terms of blocking sets, we refer to  Lemma~\ref{lemma_bs_lb_2}.

\begin{lemma}
  \label{lemma_line_mult_bound_2}
  Let $\cL$ be a multiset of lines in $\PG(n,q)$ such that each
  plane contains at most $s$ lines. Then, we have
  \begin{equation}
    \left|\cL\right| \le \left(q^4\frac{(q^{n-1}-1)(q^{n-2}-1)}{(q^3-1)(q^2-1)}+[n-1]_q\right)\cdot s-q[n-2]_q\cdot\cL(L)
  \end{equation}
  for each line $L$, where $\cL(L)$ denotes the multiplicity of $L$ in $\mathcal{L}$.
\end{lemma}
\begin{proof}
{
Fix a line $L$ and let $\cP_1$ be the set of planes that contain $L$, and  $\cP_2$ be the set of planes that are disjoint to $L$. Hence, $|\cP_1|=[n-1]_q$ and $|\cP_2|=q^6\qbin{n-1}{3}{q}$. Consider the multiset $\cP:=q^2[n-3]_q\cdot\cP_1+\cP_2$ of $q^2[n-3]_q[n-1]_q+q^6\qbin{n-1}{3}{q}$ planes. Note that $L$ is contained in all elements of $\cP_1$ and so in $q^2[n-3]_q [n-1]_q$ elements of $\cP$. Any line $L'$ that is disjoint to $L$ is contained in $q^2[n-3]_q$ elements of $\cP_2$ and $\cP$. All other lines (i.e.\ those that intersect $L$ in a point) are contained in a unique element from $\cP_1$ and so $q^2[n-3]_q$ elements from $\cP$. Consider a projective $(2,n-2)-(|\cL|,n+1,s)$ system $\cL$.

  We double count the set $S=\{(l', \pi)\,:\, l' \in \cL, l'\subset  \pi, \pi\in \cP\}$; which gives that
  \begin{align*} 
     \cL (L) q^2[n-3]_q [n-1]_q+\sum_{l' \neq L, l'\cap L\neq \emptyset} \cL(l') q^2[n-3]_q + \sum_{ l'\cap L= \emptyset} \cL(l') q^2[n-3]_q \\ \leq \left([n-1]_q q^2[n-3]_q+q^6\qbin{n-1}{3}{q}\right) s.
     \end{align*}
    This is equivalent to 
    \begin{align*}
      q^2[n-3]_q|\cL| + \cL (L)  q^3[n-3]_q[n-2]_q \leq \left([n-1]_q q^2[n-3]_q+q^6\qbin{n-1}{3}{q}\right) s.
  \end{align*} Hence, 
  \begin{align*}
      |\cL|\le s\left([n-1]_q + q^4 \frac{(q^{n-1}-1)(q^{n-2}-1)}{(q^3-1)(q^2-1)} \right) - \cL(L) q([n-2]_q),
  \end{align*}

  which proves the lemma.\qedhere
}

\end{proof}

\begin{table}[htp]
  \begin{center}
     \begin{tabular}{rrrrrr}
       \hline
       $s$ & $n_3(5,2,2;s)$ & $s$ & $n_3(5,2,2;s)$ & $s$ & $n_3(5,2,2;s)$ \\
       \hline
       1  & 28  & 6  & 465--558 & 11 & 1004--1023  \\
       2  & 105--186 & 7 & 562--651 & 12 & 1107             \\
       3  & 190--279 & 8  & 660--744 & 13 & 1210            \\
       4  & 275--372 & 9  & 810--837            \\
       5  & 366--465 & 10 & 904--930 \\          
       \hline
     \end{tabular}
     \caption{Bounds for $n_3(5,2,2;s)$.}
     \label{table_bounds_n_3_5_2_2}
   \end{center}
 \end{table}

From Theorem~\ref{thm_eisfeld_characterization}, Lemma~\ref{lemma_q_fold_blocking_set_lines_vs_planes}, Theorem~\ref{thm_n_q_5_2_2_1}, Lemma~\ref{lemma_bs_lb_1}, and Proposition~\ref{prop_subspace_blocking} we conclude $b_3(5,2,2;1,1)=103$, $b_3(5,2,2;3,1)\le 306$, $b_3(5,2,2;12,1)= 1182$, and $b_3(5,2,2;13,1)=1210$, respectively. For 
$b_3(5,2,2;2,1)\le 206$, $b_3(5,2,2;5,1)\le 550$, $b_3(5,2,2;6,1)\le 648$, $b_3(5,2,2;7,1)\le 745$, $b_3(5,2,2;8,1)\le 844$, 
$b_3(5,2,2;9,1)\le 935$, 
$b_3(5,2,2;10,1)\le 1020$, and $b_3(5,2,2;11,1)\le 1105$ we refer to Appendix~\ref{sec_bs_lists}. Using Lemma~\ref{lemma_remove_one_blocking_set} we obtain the lower bounds for $n_3(5,2,2;s)$ for $1\le s\le 13$, as summarized in Table~\ref{table_bounds_n_3_5_2_2}.   

\section{Conclusion and open problems}
\label{sec_conclusion}

We have introduced the maximum number $n_q(r,h,f;s)$ of $(h-1)$-spaces in $\PG(r-1,q)$ such that each subspace of codimension $f$ contains at most $s$ elements. These numbers are complemented by the minimum number $b_q(r,h,f;s)$ of $(h-1)$-spaces in $\PG(r-1,q)$ such that each subspace of codimension $f$ contains at least $s$ elements. Both notions are rather general. As an example, the case $(h,f)=(1,1)$ corresponds to linear codes with their geometric reformulation as multisets of points. If we keep $f=1$ but consider $h>1$, then we are dealing with additive codes. For $h=1$ and $f>1$ we are confronted with linear codes w.r.t.\ to the $f$th generalized Hamming weight. So, in this paper, we generalize both concepts to one more general structure. Due to this general setting, one cannot expect to determine these number in full generality. While we have some results on the asymptotic behavior, even the question for the right order of magnitude remains open in most cases.  Besides a few general insights we mostly focused on $n_q(5,2,2;s)$ and $b_q(5,2,2;s)$, where we mostly assume $q\in \{2,3\}$. As a first specific open problem we ask for the right order of magnitude of $n_q(5,2,2;2)$ in terms of $q$.

In Theorem~\ref{thm_b_2_5_2_2} we have fully determined the minimum number $b_2(5,2,2;s)$ of lines in $\PG(4,2)$ such that each plane contains at least $s$ elements as a function of $s$. However, this result is still based on integer linear programming computations and we propose it as an open problem to replace some of these by theoretical lower bounds. The techniques used in \cite{eisfeld1997blocking,metsch2004blocking} may serve as a blueprint. If we restrict the maximum multiplicity of the lines, then in most cases we only presented upper bounds by listing explicit examples found by ILP searches. It would be interesting to determine the exact values. For $b_3(5,2,2;s)$ we have presented partial results, see Table \ref{table_b_3_5_2_2}.

In Theorem~\ref{thm_5_2_2_q_gen} we have fully determined $n_q(5,2,2;t\cdot[3]_q+q^2+q)$. 
The underlying construction fits into the framework of Lemma~\ref{lemma_remove_blocking_sets}:  starting from the set of all lines in $\PG(4,q)$ we can remove any set of lines that blocks all planes to obtain a lower bound for $n_q(5,2,2;q^2+q)$. 
Choosing the trivial blocking set consisting of all $\qbin{4}{2}{q}=q^4+q^3+q^2+q+1$ lines in a solid yields $n_q(5,2,2;q^2+q)\ge q^6+q^5+q^4+q^3$, i.e.\ $n_2(5,2,2;6)\ge 120$. 
Choosing the blocking set obtained from the $q^4+q^3+q^2+q+1$ lines in the orbit of a Singer-cycle of $\PG(4,q)$ yields $n_q(5,2,2;q^2+q)\ge q^6+q^5+q^4+q^3+q^2$, i.e.\ $n_2(5,2,2;6)\ge 124$. 
The best choice of the blocking set yields the lower bound from Theorem~\ref{thm_5_2_2_q_gen}, i.e.\ $n_2(5,2,2;6)\ge 128$, which is tight. 
So far, all of our lower bounds for $n_2(5,2,2;s)$ are of this type. Finding a good lower bound for $n_q(5,2,2;2)$ seems to be a challenging problem.

While there is a Griesmer type bound for linear and additive codes that determines $n_q(r,h,1;s)$ for all sufficiently large values of $s$, we currently do not know such a bound for the cases $h,f\ge 2$.

In order to turn the determination of $n_q(5,2,2;s)$ and $b_q(5,2,2;s)$ as a function of $s$, given some fixed field size $q$, into a finite computational problem, we have presented Lemma~\ref{lemma_line_mult_bound_1} and Lemma~\ref{lemma_bs_lb_1}. Both bounds are generalized to some extent, but still do not cover the whole parameter space of $(r,h,f)$. We can conclude that in this paper, we give a new, rather general research direction, in which many things can still be investigated.

\section*{Acknowledgements}
The authors would like to thank Timothy Alderson, Simeon Ball, and Tabriz Popatia for the discussions on additive codes during the seventh Irsee conference. There, we uncovered 
the relation between the geometric objects we study in this paper, and additive codes with respect to the generalized Hamming weight. Both authors discussed the initial ideas 
for this paper at that conference.

\bibliography{bib}
\bibliographystyle{abbrv}

\appendix
\section{Explicit lists of blocking sets found via ILP searches}
\label{sec_bs_lists}
In this section we collect a list of interesting blocking sets that we have found by integer linear programming computations using the \texttt{ILOG CPLEX} solver without any symmetry reductions or specific settings. Some of these examples show that there are no uniqueness results for certain parameters. Others have larger cardinalities than the optimum blocking sets but require a smaller maximum line multiplicity.
For each example we state an explicit list of generator matrices of all involved subspaces.

\medskip 

A blocking set attaining cardinality $b_2(5,2,2;3)=75$ with maximum line multiplicity $3$ is given by: $\left(
\right)$.

\smallskip

A blocking set attaining cardinality $b_2(5,2,2;3)=75$ with maximum line multiplicity $2$ is given by:
$\left(%
\right)$. There are $21$ lines of multiplicity $2$ and $56$ lines of multiplicity $1$.
\smallskip

$b_2(5,2,2;5,1)\le 123$:
$\left(%
\right)$. The seven lines of multiplicity $3$ are the lines contained in a plane $\pi$. The $48$ lines of multiplicity $0$ intersect $\pi$ in a point and there are no lines of multiplicity $2$.

\smallskip 
$b_2(5,2,2;5,4)\le 120$:
$\left(%
\right)$. The seven lines of multiplicity $4$ are the lines of a plane $\pi$. The eight lines of multiplicity $2$ are disjoint to $\pi$ and the $64$ lines of multiplicity $0$ intersect $\pi$ in a point.

\smallskip 

The complement of a partial line spread of size $9$ gives $b_2(5,2,2;6,1)\le 146$:
$\left(%
\right)$. Five of the seven lines of multiplicity $2$ are contained in a plane. 

\smallskip 

$b_2(5,2,2;6,3)\le 142$:
$\left(%
\right)$. There is a unique line of multiplicity $5$ and all other lines have multiplicity at most $2$.

\smallskip

A blocking set with line multiplicity $3$ showing $b_3(5,2,2;5)\le 502$ is given by:
$\left(%
\right)$.

\smallskip

A blocking set with line multiplicity $6$ showing $b_3(5,2,2;6)\le 600$ is given by:
$\left(%
\right)$.

\smallskip

A blocking set showing $b_3(5,2,2;2,1)\le 206$ is given by:
$\left(%
\right)$.

\smallskip

A blocking set showing $b_3(5,2,2;5,1)\le 550$ is given by:
$\left(%
\right)$.

\smallskip

A blocking set showing $b_3(5,2,2;6,1)\le 648$ is given by:
$\left(%
\right)$.

\smallskip

A blocking set with line multiplicity $4$ showing $b_3(5,2,2;7)\le 690$ is given by:
$\left(%
\right)$.

\smallskip

A blocking set with line multiplicity $4$ showing $b_3(5,2,2;8)\le 784$ is given by:
$\left(%
\right)$.

\end{document}